\documentclass[a4paper,11pt]{article}

\usepackage{times,amsmath,amssymb,amsfonts}
\usepackage{graphics}
\usepackage[latin1]{inputenc}

\newtheorem{theorem}{Theorem}[section]
\newtheorem{proposition}[theorem]{Proposition}

\newtheorem{definition}[theorem]{Definition}

\rm
\rm
\rm

\begin{document}

\title{Canonical lossless state-space systems: staircase forms
and the Schur algorithm}

\author{\begin{tabular}{cccc}
Ralf L.M. Peeters & Bernard Hanzon & Martine Olivi&\tabularnewline
Dept. Mathematics & School of Mathematical Sciences & Projet
APICS&\tabularnewline Universiteit Maastricht & University College
Cork & INRIA, Sophia Antipolis&\tabularnewline
ralf.peeters@math.unimaas.nl & b.hanzon@ucc.ie &
martine.olivi@sophia.inria.fr &
\end{tabular} }

\date{}

\maketitle \thispagestyle{empty}

\begin{abstract}
A new finite atlas of overlapping balanced canonical forms for
multivariate discrete-time lossless systems is presented. The
canonical forms have the property that the controllability matrix
is positive upper triangular up to a suitable permutation of its
columns. This is a generalization of a similar balanced canonical
form for continuous-time lossless systems. It is shown that this
atlas is in fact a finite sub-atlas of the infinite atlas of
overlapping balanced canonical forms for lossless systems that is
associated with the tangential Schur algorithm; such canonical
forms satisfy certain interpolation conditions on a corresponding
sequence of lossless transfer matrices. The connection between
these balanced canonical forms for lossless systems and the
tangential Schur algorithm for lossless systems is a
generalization of the same connection in the SISO case that was
noted before. The results are directly applicable to obtain a
finite sub-atlas of multivariate input-normal canonical forms for
stable linear systems of given fixed order, which is minimal in
the sense that no chart can be left out of the atlas without
losing the property that the atlas covers the manifold.
\end{abstract}

\paragraph{Keywords:} Lossless systems, input normal forms, output normal forms,
balanced canonical forms, model reduction, MIMO systems,
tangential Schur algorithm.

\bigskip
\noindent

\section{Introduction}

In linear systems theory there has been a longstanding program in
developing balanced realizations, balanced canonical forms and
associated parameterizations for stable linear systems and for
various other classes of linear systems. The classical Gramian
based concept of balancing as introduced by Moore, see
\cite{Moore}, applies to \emph{stable} systems and allows one to
develop parameterizations in which system stability is a built-in
property. One of the motivations for the interest in balancing is
that it leads to a simple method for model order reduction, namely
by truncation of (the last entries of) the state vector.

However truncation does not always lead to a minimal system.
Therefore there has been research into balanced canonical forms
which do have the property that truncation of the last entries in
the state vector leads to a minimal system. For
\emph{continuous-time} systems this has led to the original
balanced canonical form of Ober (see \cite{Obe87}) and to the new
balanced canonical form of Hanzon (see \cite{Hanzon1995}; see also
\cite{Ober1996}). This last balanced canonical form is based on
the idea that if the controllability matrix is positive upper
triangular (i.e., the controllability matrix forms an upper
triangular matrix with positive entries on the pivot positions),
then truncation of the last entries of the state vector leads
again to a system with a positive upper triangular controllability
matrix, hence is controllable. Because this is in the balanced
continuous-time case, the controllability property here implies
that the resulting system is again minimal and balanced.

To use similar ideas to build overlapping balanced canonical forms
is more involved. For continuous-time {\em lossless} systems,
which form the key to these problems, a generalization of positive
upper triangular matrices is used in \cite{H-O98}. The idea used
there is that it suffices if a column permutation of the
controllability matrix is positive upper triangular. Under certain
circumstances there will exist an associated column permutation
(we also speak of a {\it shuffle} of columns in this context) of
the so-called realization matrix, which allows one to proceed with
the construction.

In the case of \emph{discrete-time} systems the situation is
somewhat more complicated because it is known that starting from a
balanced realization, truncation of the state vector will normally
not lead to a balanced state-space system. In the case of SISO
lossless discrete-time systems a balanced canonical form with a
simple positive upper triangular controllability matrix was
presented in \cite{HP2000}. Also the possibilities for model
reduction by truncation, combined with a correction of some sort
to arrive at a balanced realization of a lossless system, are
discussed there.

In the current paper we treat the case of MIMO lossless
discrete-time systems. We present overlapping balanced canonical
forms which have the property that the corresponding
controllability matrix is positive upper triangular, up to a
column permutation. In this sense it is close to the corresponding
results in \cite{H-O98}; however, here a generalization is
presented which simplifies the presentation and which can, as a
spin-off, also be used in the continuous-time case. The precise
relation with the approach taken in \cite{H-O98} will be made
clear. The results on the relation between a specific triangular
pivot structure in controllable pairs, which we call ``staircase
forms'', and an associated triangular pivot structure in the
controllability matrices are also of interest outside the context
of lossless systems.

In \cite{HP2000} a connection was shown between the balanced
canonical forms there presented and the Schur algorithm for scalar
lossless discrete-time transfer functions. In
\cite{HOP2006} it is shown how the parameterizations for
multivariable rational lossless transfer matrices by Schur
parameters, based on the so-called tangential Schur algorithm, can
likewise be lifted into parameterizations by Schur parameters of
balanced state-space canonical forms of lossless systems. One of
the main results of the current paper is to show how the atlas of
overlapping balanced canonical forms presented in this paper can
be obtained as a finite sub-atlas of the infinite atlas of
overlapping balanced canonical forms corresponding to the
tangential Schur algorithm. In fact, a certain well-specified
choice of so-called direction vectors in the tangential Schur
algorithm leads to the balanced canonical forms presented here.

Although a generalization of the results of this paper to the case
of complex-valued systems is straightforward, we shall restrict
the discussion to the case of real-valued systems only for ease of
presentation.

\section{Preliminaries}

\subsection{State space systems and realization theory}

Consider a linear time-invariant state-space system in discrete
time with $m$ inputs and $m$ outputs:
\begin{eqnarray}
  & & x_{t+1} = A x_t + B u_t, \\
  & & y_t     = C x_t + D u_t,
\end{eqnarray}
with $t \in {\mathbb Z}$, $x_t \in {\mathbb R}^n$ for some
nonnegative integer $n$ (the state space dimension), $u_t \in
{\mathbb R}^m$ and $y_t \in {\mathbb R}^m$. The matrices $A$, $B$,
$C$ and $D$ with real-valued entries are of compatible sizes: $n
\times n$, $n \times m$, $m \times n$ and $m \times m$,
respectively. The corresponding transfer matrix of this system is
given by $G(z)=D+C(zI_n-A)^{-1}B$, which is an $m \times m$ matrix
with rational functions as its entries. The controllability matrix
$K$ and the observability matrix $O$ associated with this system
are defined as the block-partitioned matrices
\begin{equation}
  K = [B, AB, \ldots, A^{n-1}B], \hspace{1cm}
  O = \left[ \begin{array}{c} C \\ CA \\ \vdots \\ CA^{n-1} \end{array} \right].
\end{equation}
The system (or its input pair $(A,B)$) is called controllable if
$K$ has full row rank $n$ and the system (or its output pair
$(C,A)$) is called observable if $O$ has full column rank $n$.
Minimality holds iff both controllability and observability hold,
which holds iff the McMillan degree of $G(z)$ is equal to $n$.

To any such state-space system we associate the following (square)
block-partitioned matrix $R$, which we call the {\em realization
matrix}:
\begin{equation}
  R = \left[ \begin{array}{cc} D & C \\ B & A \end{array} \right].
\end{equation}
The matrix $R$, its $n \times (m+n)$ sub-matrix $[B, A]$, and the
associated $n \times nm$ controllability matrix $K$ will all play
an important role in the sequel.

\subsection{Stability and balancing}

Let $(A,B,C,D)$ be some state space realization of a transfer
matrix $G(z)$. If the eigenvalues of $A$ all belong to the open
unit disk in the complex plane, then the matrix $A$ is called
(discrete-time) \emph{asymptotically stable}, and $(A,B,C,D)$ is
an asymptotically stable realization of $G(z)$. (For more details
on state-space realization theory, see e.g. \cite{Kai80}.)

If $(A,B,C,D)$ is an asymptotically stable realization, then the
controllability Gramian $W_c$ and the observability Gramian $W_o$
are well defined as the exponentially convergent series
\begin{eqnarray}
\label{Gramians}%
  & & W_c = \sum_{k=0}^{\infty} A^k B B^T (A^T)^k, \\
  & & W_o = \sum_{k=0}^{\infty} (A^T)^k C^T C A^{k}.
\end{eqnarray}
These Gramians are characterized as the unique (and positive
semi-definite) solutions of the respective Lyapunov-Stein
equations
\begin{eqnarray}
\label{Steineq1}%
  W_c - A W_c A^T & = & B B^T, \\
\label{Steineq2}%
  W_o - A^T W_o A & = & C^T C.
\end{eqnarray}
A minimal and asymptotically stable state-space realization
$(A,B,C,D)$ of a transfer matrix is called \emph{balanced} if its
controllability and observability Gramians $W_c$ and $W_o$ are
both diagonal and equal. Minimality implies that $W_c$ and $W_o$
are non-singular, hence positive definite. Any minimal and
asymptotically stable realization $(A,B,C,D)$ is similar to a
balanced realization, meaning that there exists a nonsingular
state space transformation matrix $T$ which makes the realization
$(TAT^{-1}, TB, CT^{-1},D)$ into a balanced realization.

A system is called \emph{input-normal} if $W_c=I_n$ and it is
called \emph{output-normal} if $W_o=I_n$. Balanced realizations
are directly related to input-normal and output-normal
realizations, respectively, by diagonal state space
transformations. The property of input-normality (resp.
output-normality) is preserved under orthogonal state space
transformations.

\subsection{Lossless systems, balanced realizations and the tangential
Schur algorithm}

A discrete-time system is called \emph{lossless} if it is stable
and its $m \times m$ transfer matrix $G(z)$ is unitary for all
complex $z$ with $|z|=1$. It is well-known (cf., e.g., Proposition
3.2 in \cite{HOP2006} and the references given there) that
$R = \left[ \begin{array}{cc} D & C \\ B & A \end{array} \right]$
is a balanced realization matrix of a lossless system if and only
if $R$ is an orthogonal matrix and $A$ is asymptotically stable.
It then holds that $W_c = W_o = I_n$. For a further background on
lossless systems, see e.g. \cite{Getal}.

In \cite{HOP2006} an atlas of overlapping balanced canonical
forms for lossless discrete-time systems of order $n$ is
presented. Also, a closely related atlas is given for
(controllable) input-normal pairs $(A,B)$ by considering the
quotient space with respect to the orthogonal group. Each of these
balanced canonical forms is then characterized (in the \emph{real}
case) by a fixed sequence of $n$ \emph{interpolation points} $w_k
\in {\mathbb R}$, $|w_k|<1$, $k=1,\ldots,n$, and a fixed sequence
of $n$ normalized \emph{direction vectors} $u_k \in {\mathbb
R}^m$, $\|u_k\|=1$, $k=1,\ldots,n$ (which are not to be confused
with the input signal applied to a system). Here we shall consider
the case $w_k=0$, $k=1,\ldots,n$, hence each balanced canonical
form that we consider is determined entirely by the choice of
direction vectors. Each such balanced canonical form for
input-normal pairs $(A,B)$ is then parameterized by a sequence of
$n$ \emph{Schur vectors} $v_k \in {\mathbb R}^m$, with $\|v_k\|<1$
for all $k=1,\ldots,n$. For lossless systems the parameterization
also involves an additional $m \times m$ orthogonal matrix $D_0$.

In fact, the realization matrix $R$ in this set-up can be written
as an \emph{orthogonal matrix product}:
\begin{equation}
\label{matrixproductformula}%
  R = \Gamma_n \Gamma_{n-1} \cdots \Gamma_1 R_0
      \Delta_1^T \Delta_2^T \cdots \Delta_n^T,
\end{equation}
where for $k=1,\ldots,n$:
\[
  \Gamma_k = \left[ \begin{array}{ccc} I_{n-k} & 0     & 0  \\
                                       0       & V_k   & 0  \\
                                       0       & 0     & I_{k-1}
                    \end{array} \right],
\]
\[
  \Delta_k = \left[ \begin{array}{ccc} I_{n-k} & 0     & 0  \\
                                       0       & U_k   & 0  \\
                                       0       & 0     & I_{k-1}
                    \end{array} \right]
\]
with an $(m+1) \times (m+1)$ orthogonal matrix block $V_k$ given
by
\[
  V_k = \left[ \begin{array}{cc}
  v_k & I_m-(1-\sqrt{1-\|v_k\|^2)}\frac{v_k v_k^T}{\|v_k\|^2} \\
  \sqrt{1-\|v_k\|^2} & -v_k^T \end{array} \right],
\]
(for $v_k=0$ it holds that $V_k=\left[ \begin{array}{cc} 0 & I_m \\
1 & 0 \end{array} \right]$ which makes that $V_k$ depends smoothly
on the entries of $v_k$) and an $(m+1) \times (m+1)$ orthogonal
matrix block $U_k$ given by
\[
  U_k = \left[ \begin{array}{cc}
        u_k & I_m-u_k u_k^T  \\
        0   & u_k^T
        \end{array} \right]
\]
and furthermore an $(n+m) \times (n+m)$ orthogonal matrix $R_0$
given by
\[
  R_0 = \left[ \begin{array}{cc}
        I_n & 0 \\
        0   & D_0
        \end{array} \right]
\]
in which $D_0$ is $m \times m$ orthogonal.

The interpolation conditions attain the form
$G_k(w_k^{-1})u_k=v_k$, where $G_k(z)$ denotes the transfer
function associated with the $k$-th order lossless system for
which the right lower $(m+k) \times (m+k)$ sub-matrix of $R_k =
\Gamma_k \Gamma_{k-1} \cdots \Gamma_1 R_0 \Delta_1^T \cdots
\Delta_{k-1}^T \Delta_k^T$ is a realization matrix. In the present
situation with $w_k=0$ it follows that $G_k(w_k^{-1}) = G_k(\infty)
= D_k$, so that the interpolation conditions can be written as
\[
  D_k u_k = v_k
\]
where $(A_k,B_k,C_k,D_k)$ denotes the corresponding state-space
realization of the $k$-th order lossless function $G_k(z)$.

Note that here we consider \emph{the real case} with real
direction vectors and real Schur parameter vectors. Note further
that $R_0$, $\Gamma_1,\ldots,\Gamma_n$ and
$\Delta_1,\ldots,\Delta_n$ are all orthogonal matrices. It is
important to note that the orthogonal matrix product $\Gamma_n
\Gamma_{n-1} \cdots \Gamma_1 R_0$ in fact forms a \emph{positive
$m$-upper Hessenberg matrix}, i.e. an $(m+n) \times (m+n)$ matrix
of which the $m$-th sub-diagonal has positive entries only and of
which the last $n-1$ sub-diagonals are all zero. It also follows
that if the direction vectors $u_1,\ldots,u_n$ are taken to be
standard basis vectors, then the matrix product $\Delta_1^T
\Delta_2^T \cdots \Delta_n^T$ yields a permutation matrix. Hence
in that case the balanced realization matrix $R$ is obtained as a
column permutation of an orthogonal positive $m$-upper Hessenberg
matrix.

\section{Triangular structures in controllable pairs and their controllability matrices}
\label{triangularstructures}%

It is not difficult to see that if the realization matrix $R$ is
positive $m$-upper Hessenberg, then (i) the first $n$ columns of
the partitioned $n \times (m+n)$ matrix $[B,A]$ form a positive
upper triangular matrix, i.e. an upper triangular matrix with only
positive entries on the main diagonal, and (ii) the first $n$
columns of the corresponding controllability matrix
$K=[B,AB,\ldots,A^{n-1}B]$ also form a positive upper triangular
matrix. (A matrix with this property is called a \emph{simple}
positive upper triangular matrix.) Therefore the realization is
controllable. In the discrete-time lossless case, if $R$ is
orthogonal, controllability implies that $A$ is asymptotically
stable which in turn implies that the realization is minimal.

A balanced realization of a lossless system is determined up to an
arbitrary orthogonal change of basis of the state space. The
effect of such a change of basis on the controllability matrix is
that it is pre-multiplied with an orthogonal matrix. Now it is
well-known that any nonsingular square matrix can be written as
the product of an orthogonal matrix and a positive upper
triangular matrix in a unique way (in numerical linear algebra
this is known as the QR-decomposition). If the first $n$ columns
of the controllability matrix are linearly independent then a
unique orthogonal state-space isomorphism exists which transforms
the first $n$ columns of the controllability matrix into a
positive upper triangular matrix. This determines a unique local
balanced canonical form for lossless systems. In the SISO case it
is in fact a global balanced canonical form and it is presented
and investigated in \cite{HP2000}.

In the MIMO case, the canonical form does not apply to systems
which have a non-generic Kronecker structure. This is why this is
a \emph{local} canonical form. In order to see how the concept of
requiring the first $n$ columns of the controllability matrix $K$
to be positive upper triangular can be generalized to obtain an
atlas of local canonical forms in the MIMO case, we will consider
the relation between triangular structures in the partitioned
matrix $[B,A]$ and triangular structures in the corresponding
controllability matrix $K=[B,AB,\ldots,A^{n-1}B]$. The following
definitions will turn out to be useful.

\begin{definition}
Let $n$ be a fixed positive integer. Consider a vector $v \in
{\mathbb R}^n$.
\\
(a) The vector $v$ is called a \emph{pivot vector} with a pivot at
position $k$, or a pivot-$k$ vector for short, if $ \in
\{1,\ldots,n\}$ is an integer for which the entry $v(k)$ is
strictly positive and the entries $v(j)$ with $j>k$ are all zero.
\\
(b) The vector $v$ is called a \emph{positive pivot-$k$ vector} if
it is a pivot-$k$ vector for which in addition the entries $v(j)$
with $j<k$ are all strictly positive too.
\end{definition}

\begin{definition}
For given positive integers $n$ and $r$, consider a mapping $J:
\{1,\ldots,n\} \rightarrow \{0,1,\ldots,r\}$ which is written in
short-hand notation as $J=\{j_1,j_2,\ldots,j_n\}$.
\\
(a) Associated with $J$, the mapping $J^+ : {\cal D}_J^+
\rightarrow {\cal R}_J^+$ is defined as the restriction of $J$ to
${\cal D}_J^+$ which is the largest subset of $\{1,\ldots,n\}$ on
which $J$ is nonzero; the co-domain ${\cal R}_J^+$ is the
corresponding range of positive values occurring as images under
$J$.
\\
(b) The mapping $J$ is called a \emph{pivot structure} if $J^+$ is
a bijection. Then the inverse of $J^+$ is denoted by $Q^+$ and the
extended mapping $Q : \{1,\ldots,r\} \rightarrow \{0,1,\ldots,n\}$
is written in short-hand notation as $Q=\{q_1,q_2,\ldots,q_r\}$
and defined by: $q_k=Q^+(k)$ for $k \in {\cal R}_J^+$, and $q_k=0$
otherwise.
\\
(c) An $n \times r$ matrix $M$ is said to have a pivot structure
$J$ if for each $k \in {\cal D}_J^+$ it holds that column $j_k$ of
$M$ is a pivot-$k$ vector. (Equivalently, each column $\ell$ of
$M$ is a pivot-$q_{\ell}$ vector, where `a pivot-$0$ vector' is
synonymous to `not a pivot vector'.)
\\
(d) A pivot structure $J$ is called a \emph{full pivot structure}
if ${\cal D}_J^+ = \{1,\ldots,n\}$.
\end{definition}
{\bf Example.~~} Let $n=5$ and $r=8$. Consider the mapping $J:
\{1,\ldots,5\} \rightarrow \{0,1,\ldots,8\}$ given by $J =
\{j_1,j_2,j_3,j_4,j_5\} = \{7,1,5,3,6\}$. It follows that the
domain and co-domain of $J^+$ are given by ${\cal D}_J^+ =
\{1,2,3,4,5\}$ and ${\cal R}_J^+ = \{1,3,5,6,7\}$, respectively.
Note that $J^+$ is a bijection, so that $J$ is a pivot structure.
Since ${\cal D}_J^+ = \{1,\ldots,n\}$ it holds that $J$ defines a
\emph{full} pivot structure. The mapping $Q$, which extends the
inverse mapping $Q^+$ of $J^+$, is given by: $Q =
\{q_1,q_2,q_3,q_4,q_5,q_6,q_7,q_8\} = \{2,0,4,0,3,5,1,0\}$. Any $5
\times 8$ matrix $M$ which has the full pivot structure $J$ is of
the following form:
\[
  M = \left[ \begin{array}{cccccccc}
  \ast & \ast & \ast & \ast & \ast & \ast & + & \ast \\
  +    & \ast & \ast & \ast & \ast & \ast & 0 & \ast \\
  0    & \ast & \ast & \ast & +    & \ast & 0 & \ast \\
  0    & \ast & +    & \ast & 0    & \ast & 0 & \ast \\
  0    & \ast & 0    & \ast & 0    & +    & 0 & \ast
  \end{array} \right]
\]
where the entries denoted by $\ast$ are allowed to have an
arbitrary value and the entries denoted by $+$ are required to be
(strictly) positive. Note that $J$ addresses the entries denoted
by $+$ for each row and $Q$ specifies the same entries for each
column.
\\[5mm]
The construction of $Q$ from a given pivot structure $J$ induces a
mapping $T_{n,r}: J \mapsto Q$. From the fact that $J^+$ and $Q^+$
are each others inverse, is not difficult to see that $T_{r,n}$
provides the inverse of $T_{n,r}$. The sets $\{(k,j_k) \, | \, k
\in {\cal D}_J^+ \}$ and $\{(q_{\ell},\ell) \, | \, \ell \in {\cal
R}_J^+ \}$ obviously coincide: the mappings $J$ and $Q$ both serve
to specify the \emph{same} set of pivot locations for an $n \times
r$ matrix, satisfying the rule that in each row and in each column
of that matrix at most one pivot location is selected. The mapping
$J$ specifies these pivot locations in a row-oriented fashion, the
mapping $Q$ in a column-oriented fashion.

For a full pivot structure it holds that $n \leq r$. If
$J=\{j_1,j_2,\ldots,j_n\}$ is a full pivot structure for an $n
\times r$ matrix $M$, then the ordered selection of columns
$j_1,j_2,\ldots,j_n$ from $M$ constitutes a positive upper
triangular $n \times n$ sub-matrix. In this way, positive upper
triangularity is generalized by the concept of a full pivot
structure.

As explained before, if a block-partitioned $n \times (m+n)$
matrix $[B,A]$ is simple positive upper triangular (i.e., it has
the full pivot structure $J$ with $j_k=k$ for $k=1,\ldots,n$) then
the associated (finite or infinite) controllability matrix
$K=[B,AB,A^2 B,\ldots]$ also is simple positive upper triangular.
We now proceed to investigate the question which full pivot
structures for $[B,A]$ induce full pivot structures for $K$.
Conversely, it is of interest to determine which full pivot
structures for $K$ are induced by full pivot structures for
$[B,A]$. The latter question is more involved and it is studied in
detail in the following section. Here we address the former
question for which the following definition is useful.
\begin{definition}
Let $m$ and $n$ be given positive integers.
\\
(a) A pivot structure $F$ for an $n \times n$ matrix $A$ is called
a \emph{staircase form} for $A$ if $F^+$ is monotonically
increasing having the range ${\cal R}_F^+ = \{1,2,\ldots,p_A\}$.
Here $p_A$ denotes the number of pivots, i.e. the number of
elements in ${\cal D}_F^+$.
\\
(b) A pivot structure $J=\{j_1,\ldots,j_n\}$ for an $n \times
(m+n)$ block-partitioned matrix $[B,A]$ induces a pivot structure
$P=\{p_1,\ldots,p_n\}$ for the matrix $A$ as given by
$p_k=\max\{j_k-m,0\}$ for $k=1,\ldots,n$.
\\
(c) A full pivot structure $J=\{j_1,j_2,\ldots,j_n\}$ for an $n
\times (m+n)$ block-partitioned matrix $[B,A]$ is called an
\emph{admissible pivot structure} for $[B,A]$ if it holds that:
(i) $B$ has a pivot-$1$ vector, i.e. $1 \leq j_1 \leq m$, and (ii)
the pivot structure $P$ induced by $J$ constitutes a staircase
form for $A$.
\end{definition}
Of course, a pivot structure $J=\{j_1,\ldots,j_n\}$ for an $n
\times (m+n)$ block-partitioned matrix $[B,A]$ also induces a
pivot structure for the matrix $B$. For several purposes, the
induced pivot structures for $A$ and $B$ are more conveniently
described in terms of the associated column-oriented description
$Q=\{q_1,\ldots,q_{m+n}\}$ for $[B,A]$. For the matrix $A$ it
holds that the associated column-oriented pivot structure
$S=\{s_1,\ldots,s_n\}$ satisfies $s_k=q_{m+k}$ for all
$k=1,\ldots,n$. For the matrix $B$ the associated column-oriented
pivot structure is the restriction of $Q$ to the domain
$\{1,\ldots,m\}$, simply described by the sequence
$\{q_1,\ldots,q_m\}$.
\\[5mm]
\noindent%
{\bf Example.~~} Let $m=4$, $n=6$ and consider the full pivot
structure $J=\{3,1,5,6,4,7\}$ for the $n \times (m+n)$ partitioned
matrix $[B,A]$. The corresponding column-oriented description is
given by $Q=\{2,0,1,5,3,4,6,0,0,0\}$. The matrix $[B,A]$ therefore
has the form:
\[
  [B,A] = \left[ \begin{array}{cccc|cccccc}
  \ast & \ast & + & \ast & \ast & \ast & \ast & \ast & \ast & \ast \\
  +    & \ast & 0 & \ast & \ast & \ast & \ast & \ast & \ast & \ast \\
  0    & \ast & 0 & \ast & +    & \ast & \ast & \ast & \ast & \ast \\
  0    & \ast & 0 & \ast & 0    & +    & \ast & \ast & \ast & \ast \\
  0    & \ast & 0 & +    & 0    & 0    & \ast & \ast & \ast & \ast \\
  0    & \ast & 0 & 0    & 0    & 0    & +    & \ast & \ast & \ast
  \end{array} \right]
\]
The induced pivot structure for the matrix $A$ is given by $P =
\{p_1,p_2,p_3,p_4,p_5,p_6\} = \{0,0,1,2,0,3\}$, which follows from
$p_k=\max\{j_k-4,0\}$ for $k=1,2,\ldots,6$. The associated
column-oriented description is then given by $S =
\{s_1,s_2,s_3,s_4,s_5,s_6\} = \{3,4,6,0,0,0\} =
\{q_5,q_6,q_7,q_8,q_9,q_{10}\}$. The function $P^+$ is given by
the pairs $(3,1)$, $(4,2)$ and $(6,3)$; the inverse $S^+$ is given
by $(1,3)$, $(2,4)$ and $(3,6)$. Clearly, $P^+$ is monotonically
strictly increasing (and equivalently $S^+$ is monotonically
strictly increasing) so that $P$ is a staircase form for $A$. This
is clearly illustrated by the pattern constituted by the entries
denoted by $+$ in the matrix $A$ above. Also, the matrix $B$ has a
pivot-$1$ column as its third column. Therefore, $J$ constitutes
an admissible pivot structure for $[B,A]$. Note that the
column-oriented description of the pivot structure for the matrix
$B$ follows from the restriction of $Q$ as: $\{q_1,q_2,q_3,q_4\} =
\{2,0,1,5\}$.
\\[5mm]
Note that an \emph{admissible} pivot structure $J$ for $[B,A]$ is
totally determined by the induced pivot structure for $B$. In that
case, the pivot structure $S$ for the matrix $A$ having a
staircase form is given by $\{s_1,\ldots,s_n\} =
\{q_{m+1},\ldots,q_{m+p_A},0,\ldots,0\}$ where the subsequence
$\{q_{m+1},\ldots,q_{m+p_A}\}$ is positive and monotonically
increasing, consisting of the elements of $\{1,2,\ldots,n\}$ not
occurring in $\{q_1,\ldots,q_m\}$. For admissibility, the only
condition on the column-oriented pivot structure
$\{q_1,\ldots,q_m\}$ for $B$ is that $1$ occurs in this sequence.

If $v$ is a pivot-$k$ vector and $J$ is admissible, then the
staircase structure of $A$ implies that $w=Av$ is a pivot-$s_k$
vector. For this reason, the function $S$ will be called the
\emph{successor} function. (For convenience we also define
$S(0)=0$ and we recall that the terminology `a pivot-$0$ vector'
is synonymous to `not a pivot vector'.) The sequence of pivot
positions for the vectors $v, A v, A^2 v, A^3 v, \ldots$ is then
given by $k, S(k), S^2(k), S^3(k), \ldots$. Conversely, the
induced pivot structure $P$ for $A$ is called the
\emph{predecessor} function (here we also introduce $P(0)=0$).
Recall that $S^+$ and $P^+$ are each others inverse.

We have the following result.
\begin{theorem}
\label{firstpivotsresult}%
Let $m$ and $n$ be given positive integers.
\\
(a) If $J$ is an admissible pivot structure for an $n \times
(m+n)$ block-partitioned matrix $[B,A]$, then $K=[B,AB,A^2
B,\ldots]$ has a full pivot structure $\widetilde{J}$.
\\
(b) For every non-admissible full pivot structure $J$ there exists
an $n \times (m+n)$ matrix $[B,A]$ having the full pivot structure
$J$, for which $K=[B,AB,A^2 B,\ldots]$ does not have a full pivot
structure.
\end{theorem}
\noindent%
{\bf Proof.~~} (a) Admissibility of $J$ implies that $B$ has a
pivot-$1$ column. Thus, the (infinite) controllability matrix $K$
also has a pivot-$1$ column, because $B$ is a sub-matrix of $K$.
Now consider the induction hypothesis that the controllability
matrix is known to have pivots at positions $1,2,\ldots,k$, with
$1 \leq k < n$. From the admissible pivot structure of $[B,A]$,
either $B$ or $A$ has a pivot-$(k+1)$ column, depending on the
value of $j_{k+1}$. If $j_{k+1} \leq m$, then this column is in $B$
hence it also appears in $K$. Otherwise, column $j_{k+1}$ of
$[B,A]$ is in fact column $p_{k+1}=j_{k+1}-m$ of $A$. Equivalently,
$s_{\ell}=q_{m+\ell}=k+1$ for $\ell=p_{k+1}$. Because of the
staircase structure of $A$ (and because the prescribed pivot-$1$
column is in $B$) it holds that $\ell \leq k$. Since $K$ has a
pivot-$\ell$ column according to the induction hypothesis, the
matrix product $AK$ now has a pivot-$(k+1)$ column because of the
staircase structure of $A$. But $AK$ is a sub-matrix of $K$,
whence it follows that $K$ has a pivot-$(k+1)$ column. This shows
the induction step. Hence the controllability matrix has a full
pivot structure.
\\
(b) See Appendix \ref{App:proof1}.
$\mbox{}$ \hfill $\Box$
\\[5mm]
{\bf Remarks.~}%
\\
(i) For an admissible pivot structure $J$ for $[B,A]$ there is a
uniquely determined full pivot structure $\widetilde{J}$ which
applies to every controllability matrix $K$ that may occur for
each arbitrary matrix $[B,A]$ having the structure $J$. One can
easily calculate $\widetilde{J}$ using the numbered Young diagram
technique described in the following section. It is most clearly
displayed in $K$ for the example where each pivot-$k$ column in
$[B,A]$ is set equal to $e_k$ and each non-pivot column is set to
zero. \\
(ii) For given $m$ and $n$, the total number of different
admissible full pivot structures can be computed from the fact
that an admissible pivot structure is completely determined by the
pivot structure for $B$. It is given by:
$\sum_{\ell=1}^{\min\{m,n\}} \ell ! \left( \begin{array}{c} m \\
\ell \end{array} \right) \left( \begin{array}{c} n-1 \\ \ell-1
\end{array} \right)$.

\section{The Young diagram associated with an admissible pivot structure}
\label{youngdiagrams}%

Starting from an admissible pivot structure $J$ for $[B,A]$ we now
want to analyze the full pivot structure $\widetilde{J}$ for the
(finite) controllability matrix $K=[B,AB,\ldots,A^{n-1}B]$ induced
by $J$ and describe their relation.

Admissibility of $J$ implies that $1 \leq j_1 \leq m$, so that
$s_1$ is either zero (in which case $A$ has no pivots) or $s_1>1$.
Together with the staircase form of $A$ this means that for all
$k=1,\ldots,n$ either $s_k>k$ or $s_k=0$. The sequence $\{k, S(k),
S^2(k), S^3(k),\ldots \}$ therefore is strictly monotonically
increasing until at some point the value $0$ occurs after which
the sequence remains zero. This happens when $S^t(k)$ attains a
value in $\{p_A+1,\ldots,n\}$. Conversely, starting from a value
$\ell>0$ the sequence $\{\ell, P(\ell), P^2(\ell), P^3(\ell),
\ldots \}$ is strictly monotonically decreasing until at some
point the value $0$ occurs after which the sequence remains zero.
This happens when $P^t(\ell)$ attains a value in
$\{q_1,\ldots,q_m\}$.

In this way, an admissible pivot structure $J$ for $[B,A]$
generates a uniquely specified full pivot structure
$\widetilde{J}$ for the controllability matrix $K$. To visualize
this, it is helpful to introduce an $m \times n$ array
$Y=(y_{i,j})$, defined as follows: entry $y_{i,j}$ denotes the
pivot position of vector $i$ in the $j$-th block $A^{j-1}B$ of $K$
(so that $\widetilde{J}(k)=(j-1)m+i$ where $i$ and $j$ are such
that $y_{i,j}=k$). In terms of the column-oriented description
$\widetilde{Q}=\{\widetilde{q}_1, \widetilde{q}_2, \ldots,
\widetilde{q}_{nm}\}$ of the pivot structure of $K$ associated
with the row-oriented full pivot structure $\widetilde{J}$, it
simply holds that $y_{i,j}=\widetilde{q}_{(j-1)m+i}$ for all
$i=1,\ldots,m$ and $j=1,\ldots,n$. The array $Y$ can therefore be
regarded as an $m \times n$ matrix representation of
$\widetilde{Q}$ which allows a clearer expression of the role
played by the block-partitioning of $K$. Obviously, there is a
one-to-one correspondence between such an array $Y$ (with entries
in $\{0,1,\ldots,n\}$) and the function $\widetilde{Q}$ (from
$\{1,\ldots,nm\}$ to $\{0,1,\ldots,n\}$).
\\[5mm]
{\bf Example.~~} Let $m=4$, $n=6$ and consider the admissible full
pivot structure $J=\{3,1,5,6,4,7\}$ and its associated
column-oriented description $Q=\{2,0,1,5,3,4,6,0,0,0\}$ for the $6
\times 10$ partitioned matrix $[B,A]$ given by:
\[
  [B,A] = \left[ \begin{array}{cccc|cccccc}
  \ast & \ast & + & \ast & \ast & \ast & \ast & \ast & \ast & \ast \\
  +    & \ast & 0 & \ast & \ast & \ast & \ast & \ast & \ast & \ast \\
  0    & \ast & 0 & \ast & +    & \ast & \ast & \ast & \ast & \ast \\
  0    & \ast & 0 & \ast & 0    & +    & \ast & \ast & \ast & \ast \\
  0    & \ast & 0 & +    & 0    & 0    & \ast & \ast & \ast & \ast \\
  0    & \ast & 0 & 0    & 0    & 0    & +    & \ast & \ast & \ast
  \end{array} \right]
\]
Then the successor function $S$ is given by $S(0)=0$ and
$\{s_1,s_2,s_3,s_4,s_5,s_6\} = \{3,4,6,0,0,0\}$ and the
predecessor function $P$ is given by $P(0)=0$ and
$\{p_1,p_2,p_3,p_4,p_5,p_6\} = \{0,0,1,2,0,3\}$. Note that the
matrix $K=[B,AB,A^2B,\ldots]$ is of the form:
\[
  K = \left[ \begin{array}{cccc|cccc|cccc|cc}
  \ast & \ast & + & \ast & \ast & \ast & \ast & \ast & \ast & \ast & \ast & \ast & \ldots & \ldots \\
  +    & \ast & 0 & \ast & \ast & \ast & \ast & \ast & \ast & \ast & \ast & \ast & \ldots & \ldots \\
  0    & \ast & 0 & \ast & \ast & \ast & +    & \ast & \ast & \ast & \ast & \ast & \ldots & \ldots \\
  0    & \ast & 0 & \ast & +    & \ast & 0    & \ast & \ast & \ast & \ast & \ast & \ldots & \ldots \\
  0    & \ast & 0 & +    & 0    & \ast & 0    & \ast & \ast & \ast & \ast & \ast & \ldots & \ldots \\
  0    & \ast & 0 & 0    & 0    & \ast & 0    & \ast & \ast & \ast & +    & \ast & \ldots & \ldots
  \end{array} \right]
\]
This shows that the induced full pivot structure $\widetilde{J}$
for $K$ is given by $\widetilde{J} = \{3,1,7,5,4,11\}$ and it has
an associated column-oriented description $\widetilde{Q} =
\{2,0,1,5,4,0,3,0,0,0,6,0,\ldots\}$. The corresponding $4 \times
6$ array $Y$ is filled with the values in $\widetilde{Q}$ column
after column, yielding the diagram:
\[
  Y = \begin{array}{|c|c|c|c|c|c|}
  \hline
  2 & 4 & 0 & 0 & 0 & 0 \\ \hline
  0 & 0 & 0 & 0 & 0 & 0 \\ \hline
  1 & 3 & 6 & 0 & 0 & 0 \\ \hline
  5 & 0 & 0 & 0 & 0 & 0 \\ \hline
  \end{array}
\]
Note that the first column of $Y$ specifies the pivot structure of
$B$, i.e.: $\{q_1,q_2,q_3,q_4\} = \{2,0,1,5\}$. The other entries
of $Y$ satisfy the rule $y_{i,j+1}=S(y_{i,j})$.
\begin{theorem}
\label{inducedstructure}%
Let $J$ be an admissible full pivot structure for the
block-partitioned matrix $[B,A]$, with an associated
column-oriented description $Q=\{q_1,\ldots,q_{m+n}\}$ and the
successor function $S$ given by $S=\{s_1,\ldots,s_n\} =
\{q_{m+1},\ldots,q_{m+n}\}$ and $S(0)=0$. Then $J$ induces a full
pivot structure $\widetilde{J}$ for the (finite) controllability
matrix $K=[B,AB,\ldots,A^{n-1}B]$ which is specified in terms of
the $m \times n$ array $Y$ associated with $\widetilde{Q}$ as
follows: \\
(i) $y_{i,1}=q_i$ for $i=1,\ldots,m$; \\
(ii) $y_{i,j+1}=S(y_{i,j})$ for $i=1,\ldots,m$ and
$j=1,\ldots,n-1$.
\end{theorem}
\noindent%
{\bf Proof.~~} As argued in the previous section, the admissible
pivot structure $J$ for $[B,A]$ is entirely determined by the
induced column-oriented pivot structure $\{q_1,\ldots,q_m\}$ for
$B$. Given these (prescribed) pivot positions for the columns of
$B$, the resulting pivot positions for the columns of the block
$AB$ are given by $\{S(q_1),\ldots,S(q_m)\}$. Likewise, the pivot
positions for the columns of the block $A^2 B$ are given by
$\{S^2(q_1),\ldots,S^2(q_m)\}$. Proceeding in this fashion, it
follows that the pivot structure $\widetilde{J}$ for $K$ induced
by $J$ corresponds to an array $Y$ which is described by: (i) the
first column of $Y$, which corresponds to $B$ and satisfies
$y_{i,1}=q_i$ for $i=1,\ldots,m$; (ii) the other columns of $Y$,
which are given by the recursion $y_{i,j+1}=S(y_{i,j})$ for
$i=1,\ldots,m$ and $j=1,\ldots,n-1$. In part (a) of the proof of
Theorem \ref{firstpivotsresult} it has already been argued that
$\widetilde{J}$ obtained in this way describes a full pivot
structure for $K$.
$\mbox{}$ \hfill $\Box$
\\[5mm]
The array $Y$ in the theorem above has the property that the
values $1,2,\ldots,n$ all occur precisely once while the other
$(m-1)n$ entries are all zero. The set of arrays $Y$ with this
property is denoted by ${\cal Y}(m,n)$. Clearly, there is a
one-to-one correspondence between this set of arrays and the set
of full pivot structures for finite controllability matrices $K$
of size $n \times nm$. However, not all the arrays $Y$ in the set
${\cal Y}(m,n)$ are induced by some \emph{admissible} pivot
structure $J$ for $[B,A]$. The following definition serves the
goal of characterizing the subset of ${\cal Y}(m,n)$ of arrays $Y$
that are induced by admissible pivot structures.
\begin{definition}
An array $Y \in {\cal Y}(m,n)$ is called an \emph{admissible
numbered Young diagram} if it has the following three properties: \\
(i) for all $i=1,\ldots,m$ and $j=1,\ldots,n-1$ it holds that
$y_{i,j+1}>0$ implies $y_{i,j}>0$ ; \\
(ii) the values $n-p_B+1,\ldots,n$ all occur in different rows of
$Y$ as their last nonzero entries, where $p_B$ is the number of
nonzero rows of $Y$; \\
(iii) for all $i,i' = 1,\ldots,m$ and $j,j' = 1,\ldots,n-1$ it
holds that $y_{i,j+1}>y_{i',j'+1}>0$ implies
$y_{i,j}>y_{i',j'}>0$.
\end{definition}
Note that the number of nonzero rows of the array $Y$
corresponding to the induced full pivot structure $\widetilde{J}$
in Theorem \ref{inducedstructure} is equal to the number of
nonzero entries in the first column of $Y$, which is equal to the
number of pivots in the matrix $B$. This explains the notation
$p_B$ in the definition above. The terminology `numbered Young
diagram' will become more clear below, when the relationship with
nice selections and Young diagrams is explained.
\begin{theorem}
(a) Let $J$ be an admissible full pivot structure for the
block-partitioned matrix $[B,A]$, then the induced full pivot
structure $\widetilde{J}$ for the controllability matrix
$K=[B,AB,\ldots,A^{n-1}B]$ corresponds to an admissible numbered
Young diagram $Y$.
\\
(b) Let $Y$ be an admissible numbered Young diagram. Then there
exists an admissible full pivot structure $J$ for $[B,A]$ which
induces the full pivot structure $\widetilde{J}$ for
$K=[B,AB,\ldots,A^{n-1}B]$ which corresponds to $Y$.
\end{theorem}
\noindent%
{\bf Proof.~~} (a) From Theorem \ref{inducedstructure} we have
that $J$ induces the full pivot structure $\widetilde{J}$ for $K$
which corresponds to an array $Y \in {\cal Y}(m,n)$ given by: (i)
$y_{i,1}=q_i$ for $i=1,\ldots,m$; (ii) $y_{i,j+1}=S(y_{i,j})$ for
$i=1,\ldots,m$ and $j=1,\ldots,n-1$. Clearly, the $i$-th row of
$Y$ is entirely zero if and only if $q_i=0$. Hence the number of
nonzero rows of $Y$ is equal to the number $p_B$ of (prescribed)
pivots in $B$. As we have seen, admissibility of $J$ implies that
$s_k=0$ if and only if $k \in \{p_A+1,\ldots,p_A+p_B=n\}$. This
shows that the last nonzero entries in the $p_B$ nonzero rows of
$Y$ have the values $n-p_B+1,\ldots,n$ and they necessarily all
occur in different rows. Next, if $y_{i,j+1}>0$, then
$y_{i,j+1}=S(y_{i,j})$ with $y_{i,j}>0$ because $S(0)=0$. This
relationship is also described by the predecessor function $P$ as
$y_{i,j}=P(y_{i,j+1})>0$. Note that in fact the restricted
positive functions $S^+$ and $P^+$ describe this relationship and
they are both strictly monotonically increasing because of the
staircase property of $A$. Therefore, by application of $P^+$, the
relationship $y_{i,j+1}>y_{i',j'+1}>0$ implies that
$y_{i,j}>y_{i',j'}>0$. This shows that $Y$ is an admissible
numbered Young diagram.
\\
(b) Suppose that $Y \in {\cal Y}(m,n)$ is an admissible numbered
Young diagram. Consider the $p_B$ nonzero rows of $Y$. According
to property (ii), the last nonzero entries of these rows precisely
cover the range $\{n-p_B+1,\ldots,n\}$. It follows that all the
other entries of $Y$ are $\leq n-p_B$ because every positive value
from $\{1,\ldots,n\}$ occurs exactly once. Now consider the
function $S: \{0,1,\ldots,n\} \rightarrow \{0,1,\ldots,n\}$
defined from the values in $Y$ as follows: $S(0)=0$,
$S(y_{i,j})=y_{i,j+1}$ for all $i=1,\ldots,m$ and
$j=1,\ldots,n-1$, and $S(y_{i,n})=0$ for all $i=1,\ldots,m$. Note
that the pattern of positive values in the array $Y$ is
left-aligned according to property (i). This makes that the
definition $S(0)=0$ is consistent with the prescription
$S(y_{i,j})=y_{i,j+1}$ in situations where $y_{i,j}=0$, and also
with the prescription $S(y_{i,n})=0$ in situations where
$y_{i,n}=0$. Note also that $S(k)>0$ for all $k=1,\ldots,n-p_B$
and $S(k)=0$ for $k=n-p_B+1,\ldots,n$ (as well as for $k=0$). The
associated function $S^+$ is a bijection with domain
$\{1,\ldots,n-p_B\}$.

Property (iii) of $Y$ now implies that $S^+$ is monotonically
increasing. To see this, choose positive integers $k$ and $\ell$
with $S(k)>S(\ell)>0$. Then choose the unique integers $i$, $i'$,
$j$ and $j'$ such that $y_{i,j+1}=S(k)$ and $y_{i',j'+1}=S(\ell)$
and invoke property (iii) to obtain that $k>\ell>0$. Consequently,
$S^+$ can be used to prescribe a staircase form for the matrix
$A$. The $p_B$ positive values in $\{1,\ldots,n\}$ not occurring
in the range of $S^+$ are precisely those occurring in the first
column of $Y$. This first column of $Y$ serves to describe a pivot
structure for $B$. Together with $S^+$ this determines a full
pivot structure $J$ for $[B,A]$ in which $A$ has a staircase form.
For $J$ to be admissible, it remains to be shown that $B$ has a
prescribed pivot-$1$ column, or equivalently that one of the
entries in the first column of $Y$ is equal to $1$. To see this,
suppose that for some $y_{i,j}>0$ it holds that
$S(y_{i,j})=y_{i,j+1}=1$. Then the bijection $S^+$ can only be
monotonically increasing if $S(1)=0$ so that $1$ does not belong
to the domain of $S^+$, which requires $1$ to belong to the set of
$p_B$ largest values $\{n-p_B+1,\ldots,n\}$. But then $y_{i,j}>1$
also belongs to this set and occurs in a different row of $Y$,
producing a contradiction.
$\mbox{}$ \hfill $\Box$
\\[5mm]
We thus have established a bijection between admissible pivot
structures $J$ for $[B,A]$ and admissible numbered Young diagrams
$Y$ associated with $K$. To relate these results to the well-known
theory of nice selections and dynamical indices, the following
definition is useful.
\begin{definition}
Let $m$ and $n$ be given positive integers. \\
(a) The set ${\cal D}(m,n)$ is defined as the set of all
multi-indices $d=(d_1,d_2,\ldots,d_m) \in {\mathbb N}_0^m$ for
which $d_1+d_2+\ldots+d_m=n$. \\
(b) A selection of $n$ columns from an $n \times nm$
controllability matrix $K=[B,AB,\ldots,A^{n-1}B]$ is called a
\emph{nice selection} if there exists a multi-index $d \in {\cal
D}(m,n)$ for which the selected set of columns is given by $\{
A^{j-1}Be_i \, | \, j \in \{1,2,\ldots,d_i\} \mbox{~for~}
i=1,2,\ldots,m \}$.
\\
(c) A \emph{nice pivot structure} $\widetilde{J}$ for $K$ is a
full pivot structure for $K$ which constitutes a nice selection of
columns from $K$.
\\
(d) If $\widetilde{J}$ is a nice pivot structure for $K$, then the
associated multi-index $d \in {\cal D}(m,n)$ is called the vector
of dynamical indices and each number $d_i$ is called the $i$-th
dynamical index ($i=1,2,\ldots,m$) of the nice pivot structure, or
of the input pair $(A,B)$.
\end{definition}
Nice selections and vectors of dynamical indices $d$ are useful
and well-known concepts for studying the rank structures that can
be exhibited by a controllability matrix $K$. The most well-known
nice selection is the Kronecker nice selection, which consists of
the first $n$ linearly independent columns of $K$. Every nice
selection may occur as the Kronecker nice selection for some
controllability matrix $K$. (Cf., e.g., \cite{HanzonCWITracts} and
the references given there.) In the concept of nice selections
though, there are no a priori rank or linear independence
requirements and no triangularity conditions. Conversely, for a
nice pivot structure it is not required that the column selection
is a Kronecker nice selection. Note also that there are $n!$
different nice pivot structures all corresponding to the same nice
selection.

Above it has been shown that an admissible pivot structure for
$[B,A]$ induces a corresponding full pivot structure for $K$ for
which the associated array $Y \in {\cal Y}(m,n)$ is an admissible
numbered Young diagram. Conversely, all admissible numbered Young
diagrams are induced in this way. An admissible numbered Young
diagram specifies a selection of $n$ columns of $K$, which
constitutes an upper triangular sub-matrix; therefore these $n$
columns are linearly independent. From the definition of a nice
selection it should be clear that any nice selection can be
represented by an $m \times n$ binary array $Z=(z_{i,j})$ in the
following way: $z_{i,j}=1$ if column $i$ of the $j$-th block
$A^{j-1}B$ of $K$ is included in the nice selection, and
$z_{i,j}=0$ otherwise. The nonzero entries in such an array $Z$
exhibit a left-aligned pattern and the dynamical index $d_i$
denotes the number of nonzero entries in the $i$-th row of $Z$,
while $d_1+\ldots+d_m=n$. Such an array $Z$ is closely related to
the concept of a Young diagram, see \cite{bookonYoungdiagrams}. As
we have seen, any admissible numbered Young diagram $Y$ is
left-aligned and it therefore gives rise to an associated nice
selection; the induced full pivot structure $\widetilde{J}$ is a
nice pivot structure for $K$. This also explains our terminology.
For the purpose of the design of local canonical forms for various
classes of linear multivariable systems, it is important that
there exists an admissible numbered Young diagram for \emph{every}
nice selection. We therefore continue to study the relationship
between nice selections and admissible numbered Young diagrams.

Let $Z$ be a \emph{Young diagram}, i.e., a left-aligned $m \times
n$ binary array corresponding to a nice selection with an
associated vector of dynamical indices $d=(d_1,\ldots,d_m)$. A
\emph{numbered Young diagram} is obtained from $Z$ by replacing
the unit entries in $Z$ by the numbers $1,2,\ldots,n$ in some
arbitrary order, so that they all occur exactly once. The set of
$m \times n$ numbered Young diagrams is the subset of ${\cal
Y}(m,n)$ of left-aligned arrays. We will now show that for every
Young diagram $Z$ there exists an associated \emph{admissible}
numbered Young diagram $Y$. More precisely, we will characterize
\emph{all} the admissible numbered Young diagrams $Y$ that
correspond to $Z$.

To do this, it is convenient to associate with every
\emph{left-aligned} array $Y \in {\cal Y}(m,n)$ a corresponding
\emph{right-aligned} array $Y_r$ as follows. If $Y$ is
left-aligned then this means that there is an associated vector
of dynamical indices $d=(d_1,\ldots,d_m)$ such that $y_{ij}>0$
iff $j \leq d_i$. Thus, the $i$-th row of $Y$ has positive
entries at its first $d_i$ positions and zero entries at
the remaining $n-d_i$ positions. Then $Y_r$ is defined by:
$(Y_r)_{ij} = 0$ for $1 \leq j \leq n-d_i$ and $(Y_r)_{ij} =
y_{i,j-n+d_i}$ for $n-d_i+1 \leq j \leq n$. In other words: the
$d_i$ positive entries in the $i$-th row are all shifted $n-d_i$
positions to the right.
\begin{proposition}
Let $Z$ be an $m \times n$ Young diagram corresponding to a nice
selection with an associated vector of dynamical indices
$d=(d_1,\ldots,d_m)$. An $m \times n$ left-aligned array $Y$
corresponding to the same vector of dynamical indices $d$, is an
admissible numbered Young diagram if and only if there exists an
$m \times m$ permutation matrix $\Pi$ for the associated
right-aligned array $Y_r$ such that the $nm$-vector
$\mbox{vec}(\Pi Y_r) = \left( (\Pi Y_r e_1)^T, (\Pi Y_r e_2)^T,
\ldots, (\Pi Y_r e_n)^T \right)^T \in {\mathbb R}^{nm}$ obtained
by stacking the $n$ columns of the array $\Pi Y_r$, has the
property that if the zero entries are deleted then the $n$-vector
$(1,2,3,\ldots,n)^T$ is obtained.
\end{proposition}
\noindent%
{\bf Proof.~~} Suppose $Y$ is an admissible numbered Young diagram
corresponding to the vector of dynamical indices $d$. Consider the
$p_B=n-p_A$ nonzero values in the last column of $Y_r$ (where
$p_B$ denotes the number of pivots in $B$, i.e. the number of
nonzero entries in $\{q_1,\ldots,q_m\}$ which also is the number
of nonzero rows in $Y$ as well as in $Y_r$). These values
constitute a permutation of the set of values
$\{p_A+1,\ldots,p_A+p_B=n\}$. Now consider the predecessors
$\{P(p_A+1),\ldots,P(n)\}$. Note that the nonzero values among
these predecessors show up in an increasing order, because $P^+$
is monotonically increasing. Repeating the argument, it follows
that the \emph{same} permutation of the nonzero rows of $Y_r$
which makes that the nonzero entries in its last column appear in
an increasing order, achieves that such a property holds for
\emph{each} of the columns of $Y_r$. Consequently, when all the
columns of the row-permuted array $Y_r$ are stacked into a vector
with $nm$ entries using the well-known $\mbox{vec}(\cdot)$
operator, a column vector remains which is equal to $(1,2,\ldots,n)^T$
when all the zeros entries are deleted.

Conversely, starting from the given vector of dynamical indices
$d$ and an arbitrary choice of $\Pi$ permuting the nonzero rows of
$Z$, the $nm$-vector with the given property and the right-aligned
arrays $\Pi Y_r$ and $Y_r$ and the left-aligned array $Y$ are
completely determined. The left-alignment property (i) of an
admissible numbered Young diagram $Y$ is built-in. Properties (ii)
and (iii) of $Y$ are not difficult to verify either, because they
are easy for $\Pi Y_r$ and $Y_r$ and shifting the rows to move
between $Y_r$ and $Y$ does not basically change the requirements
(one only needs to take into account that zeros may occur to the
left of a string of nonzero entries in $Y_r$, but the dynamical
indices now specify the length of such a string in advance).
$\mbox{}$ \hfill $\Box$
\\[5mm]
Note that the technique used in the proof of this proposition is
constructive and can be used to generate all the admissible
numbered Young diagrams corresponding to a given nice selection.
There are $p_B!$ different possibilities, where $p_B$ can be read
off from $d$ as the number of dynamical indices $d_i>0$.
\\[5mm]
For given $m$ and $n$ and for each nice selection with a vector of
dynamical indices $d$, one can consider the family ${\cal F}(d)$ of
controllable pairs which have the additional property that the
selected columns from the controllability matrix are linearly
independent. Then we know that each controllable pair $(A,B)$ lies
in at least one of the families ${\cal F}(d)$, $d \in {\cal D}(m,n)$.

Now consider the \emph{family of all controllable pairs $(A,B)$ with
$A$ discrete-time asymptotically stable}, and the question of how to
parameterize this family up to state isomorphism. (I.e., up to
multiplication of the controllability matrix $K$ by a nonsingular
$n \times n$ matrix on the left.) Every such pair $(A,B)$ corresponds
to a positive definite controllability Gramian $W_c$, which can be
factored into $W_c=M^T M$ by making a well-defined choice for $M$,
e.g. by prescribing it to be a positive upper triangular Cholesky
factor. Such a choice can be parameterized to involve precisely
$n(n+1)/2$ independent real parameters. Using $M$ to act as a
state isomorphism transforms $(A,B)$ into an input-normal controllable
pair, but it does not affect any linear dependence relations between
the columns of the controllability matrix $K$. Hence it allows one to
restrict attention to the question of parameterizing the family of
\emph{input-normal} controllable pairs $(A,B)$ up to \emph{orthogonal}
state isomorphism.

Note that an input-normal controllable pair $(A,B)$ corresponds to a
row-orthonormal partitioned matrix $[B,A]$ for which $A$ is
asymptotically stable, and vice versa (see e.g. \cite{HOP2006},
\cite{HP2000}). Then for each admissible numbered Young diagram $Y$
the family of row-orthonormal $[B,A]$ with $A$ asymptotically stable
and with an admissible pivot structure corresponding to $Y$, forms a
local canonical form for this family. This set of local canonical
forms is covering this family in the sense that for each row-orthonormal
$[B,A]$ with $A$ asymptotically stable there exists an admissible
numbered Young diagram $Y$ and an orthogonal matrix $Q$ such that
$[QB,QAQ^T]$ has the admissible pivot structure associated with
$Y$. Furthermore, because of uniqueness of the associated
QR-decomposition due to positive upper triangularity, for such a
combination of $[B,A]$ and $Y$ the matrix $Q$ and hence
$[QB,QAQ^T]$ is unique.

An interesting question is how to obtain a minimal sub-atlas of
this atlas of local canonical forms, minimal in the sense that no
further local canonical forms can be left out without losing the
property of covering the family. To obtain a minimal
sub-atlas we have to choose one of the local canonical forms for
each $d \in {\cal D}(m,n)$. This implies that for each $d \in
{\cal D}(m,n)$ we have to choose one of the $p_B!$ possible
numberings of the associated Young diagram. As each such numbering
is associated with a permutation of the nonzero rows of the Young
diagram this choice can be fixed by specifying that permutation.
One possible choice is the unique permutation for which the
permuted dynamical indices form a non-increasing sequence, while
the order of the rows which have the same dynamical index is kept
the same. Note that this permutation is used only to determine the
\emph{numbering} in the Young diagram, the ordering of the
dynamical indices is left unchanged. With hindsight one can say
that this particular choice to obtain a minimal atlas was used in
\cite{H-O98} in a similar approach for continuous-time
input-normal pairs and lossless systems. Just as in that paper for
the continuous-time case, here each local canonical form on
discrete-time asymptotically stable input normal systems defines a
balanced local canonical form on minimal discrete-time lossless
systems of order $n$. How these balanced local canonical forms for
minimal discrete-time lossless systems of order $n$ are related to
those constructed in \cite{HOP2006} by means of the
tangential Schur algorithm is the topic of the next section.

\section{Atlases of balanced canonical forms for lossless systems}
\label{connection}%

We now have two approaches to arrive at an atlas of overlapping
balanced canonical forms for discrete-time lossless systems: one
using the balanced realizations associated with the tangential
Schur algorithm and one based on balanced realizations with an
imposed pivot structure on the row-orthonormal matrix $[B,A]$, hence
on the orthogonal realization matrix $R$. However, one of our main
results is that the second approach corresponds to making special
choices for the direction vectors in the first approach. Hence the
atlas of overlapping balanced canonical forms resulting from the
second approach is a sub-atlas of the atlas of overlapping balanced
canonical forms in the first approach. The precise formulation is
as follows.
\begin{theorem}
\label{schurshuffles}%
Let $Y$ be an admissible numbered Young diagram, corresponding to
an associated nice pivot structure $\widetilde{J}$ (for controllability
matrices) and an admissible pivot structure $J$ (for $n \times (m+n)$
matrices). For each $k=1,2,\ldots,n$, choose the direction vector
$u_{n+1-k}$ equal to $e_{i(k)}$, the $i(k)$-th standard basis vector
in ${\mathbb R}^m$, where $(i(k),j(k))$ denotes the unique pair of
indices such that $y_{i(k),j(k)}=k$. Then for any choice of the Schur
parameter vectors $v_1,v_2,\ldots,v_n$ (all of length $<1$) and
for any choice of the orthogonal matrix $D_0$, consider the $(m+n)
\times (m+n)$ orthogonal realization matrix $R$ given by
(\ref{matrixproductformula}). It follows that $J$ is an admissible
pivot structure for the sub-matrix $[B,A]$ and $\widetilde{J}$ is
a nice pivot structure for the controllability matrix $K$.
\end{theorem}
A detailed technical proof of this theorem is given in Appendix
\ref{App:proof2}.
\\[5mm]
From the point of view of the tangential Schur algorithm, it is of
interest also to directly characterize all the sequences of
direction vectors $u_1,u_2,\ldots,u_n$ that give rise to an
admissible pivot structure for the matrix $[B,A]$ (and an
accompanying nice pivot structure for the controllability matrix
$K$).
\begin{theorem}
\label{directionvectors}%
Consider a chart associated with the tangential Schur algorithm
(with all the interpolation points $w_k$ located at the origin),
specified by a sequence of direction vectors $\{u_1,\ldots,u_n\}$.
Then each $[B,A]$ resulting from this chart exhibits an admissible
pivot structure, irrespective of the choice of Schur vectors
$v_1,\ldots,v_n$, if the sequence of direction vectors consists of
standard basis vectors, say $u_k=e_{\mu(k)}$ for some indices
$\mu(1),\ldots,\mu(n)$ chosen from $\{1,2,\ldots,m\}$, satisfying
the following condition: \\
for each $k=1,2,\ldots,n-1$, if there exists a largest index
$\ell$ strictly less than $k$ such that $u_{\ell} = e_{\mu(k)}$, then
$\mu(k+1)$ is from the set $\{\mu(\ell+1),\ldots,\mu(k)\}$.
\end{theorem}
{\bf Proof.~~}%
This follows directly from the properties of the three procedures
introduced in Appendix \ref{App:proof2} to generate an admissible
numbered Young diagram. Details are left to the reader.
$\mbox{}$ \hfill $\Box$
\\[5mm]
{\bf Example.~~} Consider the same situation as for the example in
Appendix \ref{App:proof2}, where $m=5$, $n=12$ and
$(q_1,q_2,q_3,q_4,q_5)=(4,1,9,0,7)$. There it is remarked that the
choice of direction vectors $u_{n+1-k}=e_{i(k)}$ can be rewritten
as $u_k=e_{\mu(k)}$ where $\mu(k)$ denotes the index of the value
$1$ in the vector $\eta_k$, generated by the `third procedure'. In
this example it follows that the sequence
$\{u_1,u_2,\ldots,u_{12}\}$ is given by
$\{e_2,e_5,e_1,e_3,e_2,e_5,e_1,e_2,e_1,e_2,e_2,e_2\}$. Note that
this sequence satisfies the condition of Theorem
\ref{directionvectors} for all $k=1,2,\ldots,11$. E.g., for $k=6$
the previous occurrence of the vector $u_6=e_5$ happened for
$\ell=2$. The condition of the theorem requires $u_7$ to occur in
the set $\{u_3,u_4,u_5,u_6\} = \{e_1,e_3,e_2,e_5\}$, which indeed
is the case.

\section{Conclusions and discussion}

In this paper we have developed a detailed procedure to construct
an atlas of overlapping (local) balanced canonical forms for MIMO
discrete-time lossless systems $(A,B,C,D)$ of a given order $n$.
To achieve this, the concept of an admissible pivot structure for
$[B,A]$ has been introduced, which induces a full (nice) pivot
structure for the associated controllability matrix $K$. The
approach taken in Sections \ref{triangularstructures} and
\ref{youngdiagrams} has a wider range of applicability though: it
builds on the algebraic relationship between $[B,A]$ and $K$, and
it neither relies on input-normality of $(A,B)$ nor on the
discrete-time setting for the lossless systems. When one is
dealing with a system having a special (non-generic) Kronecker
structure, this can be recognized in these local canonical forms
by certain entries becoming zero. To demonstrate the structure of
the charts that constitute the atlases discussed in this paper, a
detailed worked example is included in Appendix \ref{App:atlas}.

One of the main practical reasons for studying overlapping
canonical forms is that (iterative) identification and
optimization algorithms in the MIMO case may suffer from numerical
ill-conditioning and slow convergence when they pass by systems
that not too far from systems with a non-generic structure.
Switching charts may then help to improve the algorithmic
performance. The connection in Section \ref{connection} with the
atlas of charts developed for discrete-time input-normal pairs
involving the tangential Schur algorithm is useful, because that
set-up involves (well conditioned!) orthogonal matrix computations
while it is tailored to deal with the important class of stable
systems. The tangential Schur algorithm provides one with a lot of
flexibility to design local canonical forms. Since it is
computationally expensive to switch charts at each and every
iteration, a suitably chosen finite sub-atlas is welcome. In the
present paper we have indicated the restrictions that should be
taken into account when choosing direction vectors $u_k$ from the
set of standard basis vectors, if a pivot structure is to show up
not only in $[B,A]$ but also in the controllability matrix $K$.
When a nice pivot structure is present in $K$, this has the
advantage that the impact of state vector truncation is easier to
analyze; controllability is then preserved. This is of importance
in the context of model order reduction applications.

Future research addresses the issue of monitoring the conditioning
of a chart (i.e., a local canonical form) at a given system, and
the issue of selecting a better chart when switching becomes
necessary. Since the total number of charts in an atlas quickly
grows large with the dimensions $m$ and $n$ (even for the case of
admissible pivot structures) it may not be attractive to carry out
a full search for a better chart over the entire atlas. The rank
structure in $K$ can then be instrumental in designing a quick
on-line algorithm which guarantees a certain degree of
conditioning improvement. This is currently under investigation.

\bibliographystyle{plain}
\bibliography{PeetersHanzonOlivi_LAA2007}

\appendix

\section{Proof of part (b) of Theorem \ref{firstpivotsresult}}
\label{App:proof1}%

Consider a non-admissible full pivot structure
$J=\{j_1,j_2,\ldots,j_n\}$ for $[B,A]$. Then it either holds that:
(1) $J$ does not prescribe $B$ to have have a pivot-$1$ column, or
(2) $J$ does prescribe $B$ to have a pivot-$1$ column, but $J$
does not impose a staircase structure on $A$.
\\
In case (1) it holds that $j_1>m$. We distinguish between two
situations. (i) If $j_1=m+1$, then the first column of $A$ is a
pivot-$1$ column. Then consider the following example: for each
$k=1,\ldots,n$ let column $j_k$ of $[B,A]$ be equal to the $k$-th
standard basis vector $e_k \in {\mathbb R}^n$ and let the
remaining $m$ columns of $[B,A]$ all be zero. Clearly, $[B,A]$
exhibits the given full pivot structure $J$, but $e_1^T B=0$ and
$e_1^T A =e_1^T$ so that $e_1^T K = 0$. Hence $K$ does not have a
pivot-$1$ column, so it does not have a full pivot structure.
\\
(ii) If instead $j_1>m+1$, then consider the following example:
for each $k=1,\ldots,n$ choose column $j_k$ of $[B,A]$ to be a
positive pivot-$k$ vector and let the remaining $m$ columns of
$[B,A]$ all be chosen to be strictly positive (so that effectively
they are all positive pivot-$n$ vectors). Clearly, $[B,A]$
exhibits the given full pivot structure $J$. Note that each column
in $B$ is (effectively) a positive pivot-$k$ vector with $k \geq
2$. Now, if $v$ is a positive pivot-$k$ vector, then $Av$ is a
positively weighted linear combination of the first $k$ columns of
$A$. Since all columns of $A$ are (effectively) positive
pivot-$\ell$ vectors for certain values of $\ell$, the vector $Av$
is a positive pivot-$p$ vector where $p$ is the maximal
(effective) pivot position among the first $k$ columns of $A$.
Now, the first column of $A$ has at least two nonzero entries,
because $j_1>m+1$. Therefore, each column of $AB$ is (effectively)
a positive pivot-$p$ vector with $p \geq 2$. By induction it
follows that all columns of $K$ are (effectively) positive
pivot-$p$ vectors with $p \geq 2$. Consequently, $K$ does not have
a pivot-$1$ column, so it does not have a full pivot structure.
\\
In case (2) it holds that $j_1 \leq m$, but the staircase
structure does not necessarily hold for $A$. We again distinguish
between two situations. (i) Suppose that for some $k<n$ there is a
pivot-$k$ vector in $A$ for which there is either a non-pivot
column in $A$ preceding it, or a pivot-$\ell$ vector preceding it
with $\ell>k$. Then consider basically the same example as used in
case (1) part (ii): for each $k=1,\ldots,n$ choose column $j_k$ of
$[B,A]$ to be a positive pivot-$k$ vector and let the remaining
$m$ columns of $[B,A]$ all be chosen to be strictly positive (so
that effectively they are all positive pivot-$n$ vectors). For
this example it now follows that $K$ does not have a pivot-$k$
vector, because $B$ does not have one and because for all
$p=1,\ldots,n$ the maximum (effective) pivot position among the
first $p$ columns of $A$ can never be equal to $k$.
\\
(ii) For all $k<n$ the pivot-$k$ vectors in $A$ respect the
staircase structure, but there is a prescribed pivot-$n$ vector in
$A$ which is directly preceded by a non-pivot column. If this
pivot-$n$ vector occurs in the last column of $A$, then one may
consider the same kind of example as used in case (1) part (i):
for each $k=1,\ldots,n$ let column $j_k$ of $[B,A]$ be equal to
$e_k$ and let the remaining $m$ columns of $[B,A]$ all be zero.
Now $e_n^T B=0$ and $e_n^T A =e_n^T$ so that $e_n^T K = 0$. Hence
$K$ does not have a pivot-$n$ column, so it does not have a full
pivot structure. If the pivot-$n$ vector does not occur in the
last column of $A$, then the last column of $A$ is a non-pivot
column. Summarizing, we then are in a situation where the first
$p_A-1$ columns of $A$ exhibit a staircase structure, column
$p:=j_n-m>p_A$ of $A$ is a pivot-$n$ column and the two columns
$p-1$ and $n$ of $A$ are both non-pivot columns. Then consider the
following example. Column $j_{p-1}$ of $[B,A]$ is defined as the
pivot-$(p-1)$ vector $\frac{3}{5} e_{p-1}$. Column $j_p$ of
$[B,A]$ is defined as the pivot-$p$ vector $\frac{3}{5} e_p +
\frac{16}{25} e_{p-1}$. The $3 \times 3$ sub-matrix $S$ of $A$
constituted by the intersection of its rows and columns with
indices $p-1$, $p$ and $n$ is defined as $S = \left(
\begin{array}{ccc} 0 & \frac{36}{125} & -\frac{48}{125} \\ 0 &
-\frac{12}{25} & \frac{16}{25} \\ -\frac{3}{5} & \frac{16}{25} &
\frac{12}{25} \end{array} \right)$. All remaining entries at pivot
positions in $[B,A]$ are defined to be $1$ and all other entries
are set to zero. For this example it will be shown that the
entries in the last row of $K$ are never positive, so that $K$
does not have a full pivot structure. Note that all columns in $B$
have a last entry that is equal to zero. All other columns in $K$
are of the form $Av$ for some vector $v \in {\mathbb R}^n$. For
$Av$ to have a nonzero last entry, at least one of the entries in
positions $p-1$, $p$ and $n$ of $v$ must be nonzero. Such vectors
$v$ must come from repeated pre-multiplication of the columns of
$B$ by the matrix $A$. The first vectors to have such a structure
are the pivot-$(p-1)$ vector and the pivot-$p$ vector that both
occur among the columns of $B$ and the first $p_A-1$ columns of
$A$. Once such vectors $v$ are multiplied by $A$, only the entries
in positions $p-1$, $p$ and $n$ can become nonzero: the subspace
spanned by $e_{p-1}$, $e_p$ and $e_n$ is an invariant subspace of
$A$. Restricting to this subspace, the matrix $A$ is represented
by the sub-matrix $S$ given above. Consequently, the entries in
the last row of $K$ are either zero or obtained as the entries in
the last row of the controllability matrix of the pair $(T,S)$
with $T=\left( \begin{array}{cc} \frac{3}{5} & \frac{16}{25} \\ 0
& \frac{3}{5} \\ 0 & 0 \end{array} \right)$ which represents the
pivot-$(p-1)$ and the pivot-$p$ vector in this new notation. For
the matrix $S$ it is easily established that
$S^3=\frac{544}{625}S$. It therefore suffices to compute the last
row of the matrix $[T,ST,S^2T]$, which is equal to
$(0,0,-\frac{9}{25},0,-\frac{108}{625},-\frac{36}{125})$. This
proves that all these entries are indeed non-positive.
$\mbox{}$ \hfill $\Box$%
\\[5mm]
{\bf Remark.~}
\\
In the case of balanced realizations of lossless
systems we will in addition require $[B,A]$ to have orthonormal
rows. The proof above does not entirely apply to this restricted
situation. For example if $[B,A]$ only has non-negative entries
then orthogonality of the rows requires that in each column there
is at most one nonzero entry. This requirement is violated by the
counterexamples presented in case (1) part (ii) and in case (2)
part (i) of the proof, because non-pivot columns are chosen to be
positive pivot-$n$ vectors. How to obtain a proof for this more
restricted orthonormal case is an open problem at this point. Note
however that the counterexamples presented in case (1) part (i)
and in case (2) part (ii) of the proof have in fact already been
designed to involve $[B,A]$ with orthonormal rows.

\section{Proof of Theorem \ref{schurshuffles}}
\label{App:proof2}%

To prove this theorem it is helpful first to reconsider the
precise relationship between $Y$, $\widetilde{J}$, $J$ and the
pivot structure $\{q_1,\ldots,q_m\}$ of the matrix $B$. Recall
that the number $1$ appears in the pivot structure
$\{q_1,\ldots,q_m\}$ for $B$ because of admissibility and this
sequence completely characterizes the successor function $S$, the
column-oriented description $Q$, the admissible pivot structure
$J$, the nice pivot structure $\widetilde{J}$ and the admissible
numbered Young diagram $Y$ in the way explained before.
\\[5mm]
We have previously introduced the following construction procedure
for $Y$ from $\{q_1,\ldots,q_m\}$:\\
{\bf Procedure 1}\\
(a) \emph{Construction}: the successor function $S$ is defined as
the increasing sequence of all positive numbers in
$\{1,2,\dots,n\}$ not occurring in $\{q_1,\ldots,q_m\}$ completed
by a sequence of $p_B$ zeros; in addition $S(0)=0$.\\
(b) \emph{Initialization}: set $y_{i,1}:=q_i$ for $i=1,\ldots,m$.
\\
(c) \emph{Recursion}: set $y_{i,j+1}:=S(y_{i,j})$ for
$i=1,\ldots,m$ and $j=1,\ldots,n-1$.
\\[5mm]
A second way to generate $Y$ in a dynamical fashion which avoids
the explicit construction of $S$, is by means of the following
procedure:\\
{\bf Procedure 2}\\
(a) \emph{Initialization}: set $y_{i,1}=q_i$ for $i=1,\ldots,m$.\\
(b) \emph{Recursion}: for $k=1,2,\ldots,n$, if the value $k$ has
not yet been assigned to an entry of $Y$ then select the smallest
nonzero number $y_{i,j}$ in $Y$ for which the entry $y_{i,j+1}$
immediately to its right has still not been assigned some value
and set $y_{i,j+1}=k$.\\
(c) \emph{Termination}: set all the remaining entries of $Y$ equal
to zero.\\
It is not very hard to establish that the array $Y$ constructed in
this fashion is indeed admissible and identical to the array $Y$
constructed previously with the help of $S$.
\\[5mm]
A third way to generate $Y$ from $\{q_1,\ldots,q_m\}$ employs a
sequence of vectors $\eta_k$ ($k=0,1,2,\ldots,n$) and proceeds
as follows in a backward fashion:\\
{\bf Procedure 3}\\
(a) \emph{Initialization}: set $\eta_n=(q_1,\ldots,q_m)^T$.\\
(b) \emph{Backward recursion}: for $k=n-1,\ldots,1,0$ construct
$\eta_k$ from $\eta_{k+1}$ by executing the following three rules
in the given order:\\
(1) if $(\eta_{k+1})_i=0$ then set $(\eta_k)_i:=0$;\\
(2) if $(\eta_{k+1})_i>1$ then set $(\eta_k)_i:=(\eta_{k+1})_i-1$;\\
(3) if $(\eta_{k+1})_i=1$ then define $\xi_k$ as the smallest
positive number different from all the entries of $\eta_k$ already
assigned by rules (1) and (2); if $\xi_k \leq k$ then set
$(\eta_k)_i:=\xi_k$ else set $(\eta_k)_i:=0$.\\
(c) \emph{Construction}: for each $i=1,\ldots,m$ consider the
$d_i$ values of $k$ for which $(\eta_{n+1-k})_i=1$ and assign
these values (in increasing order) to the first $d_i$ entries of
row $i$ of $Y$; set all other entries to $0$.
\\[5mm]
The validity of this third procedure for generating $Y$ can be
seen as follows. First, note that because of rule (2) in each
recursion step (b), the first column of $Y$ attains the required
form containing $q_1,\ldots,q_m$, since the number $1$ first
occurs in position $i$ of $\eta_{n+1-k}$ for $k=q_i$. Also note
that the positive integers in $\eta_n$ are all different (since
this holds for $q_1,\ldots,q_m$) and that the rules in each
recursion step (b) are such that this property is preserved for
all vectors $\eta_k$. Next, these rules are such that each vector
$\eta_k$ has precisely one entry equal to $1$, for all
$k=1,\ldots,n$. The construction in step (c) is such that all the
numbers from $\{1,\ldots,n\}$ show up precisely once in a
corresponding left-aligned numbered Young diagram $Y$. Finally,
rule (3) in each recursion step (b) guarantees that $y_{i,j}$ is
followed by $y_{i,j+1}=S(y_{i,j})$: note that $y_{i,j}=k$ is
equivalent to $(\eta_{n+1-k})_i=1$ and $(\eta_{n-k})_i=\xi_k$
implies that $y_{i,j+1}=k+\xi_k$; here $\xi_k>0$ is chosen as
small as possible, precisely in line with the second procedure for
generating $Y$.
\\[5mm]
{\bf Example.~~} Consider the situation with $m=5$, $n=12$ and
$(q_1,q_2,q_3,q_4,q_5)=(4,1,9,0,7)$. Then the successor function
$S$ is described by $\{s_1,s_2,\ldots,s_{10}\} =
\{2,3,5,6,8,10,11,12,0,0,0,0\}$. It follows that $p_B=4$ and the
corresponding admissible numbered Young diagram $Y$ is given by:
\[
  Y = \begin{array}{|c|c|c|c|c|c|c|c|c|c|c|c|}
      \hline
        4 & 6  & 10 & \phantom{12} & \phantom{12} & \phantom{12} &
        \phantom{12} & \phantom{12} & \phantom{12} & \phantom{12} &
        \phantom{12} & \phantom{12} \\ \hline
        1 & 2  &  3 & 5 & 8 & 12 &  &  &  &  &  &  \\ \hline
        9 & \phantom{12} &  &  &  &  &  &  &  &  &  &  \\ \hline
        \phantom{12} &  &  &  &  &  &  &  &  &  &  &  \\ \hline
        7 & 11 & \phantom{12} & &  &  &  &  &  &  &  &  \\ \hline
      \end{array}
\]
where all the zeros are omitted for clarity. The second procedure
for generating $Y$ without the explicit construction of the
successor function $S$ yields the same result. It proceeds from
the given first column of $Y$ by putting the value $2$ after the
value $1$, then the value $3$ after the value $2$, then the value
$5$ after the value $3$, then the value $6$ after the value $4$,
and so on. The values $1$, $4$, $7$ and $9$ are skipped, because
they have already been assigned to the first column of $Y$.

The third procedure for generating $Y$ involves the backward
recursion for the construction of the vectors $\eta_k$, for
$k=n,n-1,\ldots,1,0$. This produces the following sequence:
\[
  \begin{array}{|c|c|c|c|c|c|c|c|c|c|c|c|c|}
  \hline
  \eta_{12} & \eta_{11} & \eta_{10} & \eta_9 & \eta_8 & \eta_7 &
  \eta_6 & \eta_5 & \eta_4 & \eta_3 & \eta_2 & \eta_1 & \eta_0 \\
  \hline
        4 & 3 & 2 & 1 & 2 & 1 & 4 & 3 & 2 & 1 & 0 & 0 & 0 \\
        1 & 1 & 1 & 2 & 1 & 3 & 2 & 1 & 4 & 3 & 2 & 1 & 0 \\
        9 & 8 & 7 & 6 & 5 & 4 & 3 & 2 & 1 & 0 & 0 & 0 & 0 \\
        0 & 0 & 0 & 0 & 0 & 0 & 0 & 0 & 0 & 0 & 0 & 0 & 0 \\
        7 & 6 & 5 & 4 & 3 & 2 & 1 & 4 & 3 & 2 & 1 & 0 & 0 \\
  \hline
  \end{array}
\]
For instance, the vector $\eta_7$ is obtained from the vector
$\eta_8$ as follows. First all the entries equal to zero are
copied and all the values $(\eta_8)_i \geq 2$ are decreased by $1$
to produce the corresponding values of $(\eta_7)_i$. The value of
$(\eta_7)_2$ is addressed last, because $(\eta_8)_2=1$. At that
stage the values $1$ and $2$ have already been assigned to some
entries of $\eta_7$ and it holds that $\xi_7=3$. Because $\xi_7=3$
is not larger than the index $k=7$, this value is assigned to
$(\eta_7)_2$.

Once the vectors $\eta_{12},\eta_{11},\ldots,\eta_1$ have been
constructed, the array $Y$ is constructed by considering the
positions of the entries $1$. For the first row, these positions
are subsequently $4$, $6$ and $10$ (proceeding in the given order
from $\eta_{12}$ to $\eta_1$ corresponding to the index $n+1-k$).
For the second row we have: $1$, $2$, $3$, $5$, $8$ and $12$, and
so on.
\\[5mm]
The third way of characterizing $Y$ in terms of the pivot
structure $\{q_1,\ldots,q_m\}$ has a number of properties that are
worth noting in view of the proof of Theorem \ref{schurshuffles}
below. First, note that rule (3) in step (b) implies that the
maximum value among the entries of $\eta_k$ is at most $k$.
(Therefore, $\eta_0$ is the zero vector.) Second, all the positive
entries of a vector $\eta_k$ are different. This makes that if
$(\eta_{k+1})_i=1$ then a positive value $\xi_k$ is assigned to
$(\eta_k)_i$ for $k \geq p_B$ and the value $0$ is assigned for
$k<p_B$. Third, note that the sequence of values
$\xi_{n-1},\xi_{n-2},\ldots,\xi_{p_B}$ (in that backward order) is
increasing. Fourth, the choice of direction vectors
$u_{n+1-k}=e_{i(k)}$ can be rewritten as $u_k=e_{\mu(k)}$ where
$\mu(k)$ denotes the index of the value $1$ in the vector
$\eta_k$. Note that according to this notation,
$(\eta_k)_{\mu(k+1)}=\xi_k$ for $k=n-1,n-2,\ldots,p_B$. Finally,
it will be shown that the vectors $\eta_k$ represent the pivot
structures for the sequence of lossless systems of orders
$k=1,2,\ldots,n$ encountered in the tangential Schur algorithm for
the particular choice of direction vectors specified in Theorem
\ref{schurshuffles}.
\\[5mm]
{\bf Proof of Theorem \ref{schurshuffles}.~~} Consider the matrix
product
\[
  R = \Gamma_n \cdots \Gamma_1 R_0 \Delta_1^T \cdots \Delta_n^T.
\]
Note that the product $\Gamma_n \cdots \Gamma_1 R_0$ is positive
$m$-upper Hessenberg for any choice of Schur vectors
$v_1,\ldots,v_n$. Post-multiplication by the matrix $\Delta_1^T$
only affects the last $m+1$ columns, because the matrix
$\Delta_1^T$ is given by
\[
  \Delta_1^T = \left[ \begin{array}{ccc}
  I_{n-1} & 0 & 0 \\
  0 & e_{\mu(1)}^T & 0 \\
  0 & I_m-e_{\mu(1)}e_{\mu(1)}^T & e_{\mu(1)}
  \end{array} \right]
\]
where $\mu(1)$ denotes the location of the entry $1$ in the vector
$\eta_1$ which features in the third method for the construction
of $Y$ from $\{q_1,\ldots,q_m\}$. The precise effect is as
follows:\\
(i) column $n$ of $\Gamma_n \cdots \Gamma_1 R_0$ (having a pivot
in its last position) is moved into column $n+\mu(1)-1$;\\
(ii) column $n+\mu(1)$ of $\Gamma_n \cdots \Gamma_1 R_0$ is
moved into column $n+m$;\\
(iii) columns $n+1,\ldots,n+\mu(1)-1$ and $n+\mu(1)+1,\ldots,n+m$
of $\Gamma_n \cdots \Gamma_1 R_0$ are moved one position to the
left, into columns $n,\ldots,n+\mu(1)-2$ and
$n+\mu(1),\ldots,n+m-1$, respectively.\\
Note that the last row of $\Gamma_n \cdots \Gamma_1 R_0
\Delta_1^T$ can be regarded to have the structure:
\[
  \left[ \begin{array}{ccc} 0 & B_1 & A_1 \end{array} \right]
\]
with $A_1$ of size $1 \times 1$ and $B_1$ of size $1 \times m$.
The $1 \times (m+1)$ partitioned matrix $[B_1,A_1]$ has an
admissible pivot structure for which the column-oriented pivot
structure of $B_1$ is given by $\{0,\ldots,0,1,0,\ldots,0\}$ with
the value $1$ in position $\mu(1)$. In other words, the pivot
structure of $B_1$ is described by $\eta_1$.

Consider the last $k$ rows of the matrix product $\Gamma_n \cdots
\Gamma_1 R_0 \Delta_1^T \cdots \Delta_k^T$. Note that these can be
regarded to constitute the structure:
\[
  \left[ \begin{array}{ccc} 0 & B_k & A_k \end{array} \right]
\]
with $A_k$ of size $k \times k$ and $B_k$ of size $k \times m$.
Now suppose that the $k \times (m+k)$ partitioned matrix
$[B_k,A_k]$ is known to have an admissible pivot structure for
which the column-oriented pivot structure of $B_k$ is given by the
vector $\eta_k$. (This is the induction hypothesis.) We consider
what happens under post-multiplication by the matrix
$\Delta_{k+1}^T$. Note that this matrix is given by:
\[
  \Delta_{k+1}^T = \left[ \begin{array}{cccc}
  I_{n-k-1} & 0 & 0 & 0 \\
  0 & e_{\mu(k+1)}^T & 0 & 0 \\
  0 & I_m-e_{\mu(k+1)}e_{\mu(k+1)}^T & e_{\mu(k+1)} & 0 \\
  0 & 0 & 0 & I_k
  \end{array} \right]
\]
Therefore, post-multiplication by $\Delta_{k+1}^T$ only acts on
the columns $n-k,\ldots,n-k+m$ of the matrix $\Gamma_n \cdots
\Gamma_1 R_0 \Delta_1^T \cdots \Delta_k^T$.

The partitioned matrix $[B_{k+1},A_{k+1}]$ is then formed as
\[
  [B_{k+1},A_{k+1}] =
  \left[ \begin{array}{ccc}
    \gamma & \beta & \alpha \\
    0 & B_k & A_k
  \end{array} \right]
  \left[ \begin{array}{ccc}
    e_{\mu(k+1)}^T & 0 & 0 \\
    I_m-e_{\mu(k+1)}e_{\mu(k+1)}^T & e_{\mu(k+1)} & 0 \\
    0 & 0 & I_k
  \end{array} \right]
\]
where $\gamma$ is a positive scalar and $\beta$ and $\alpha$ are
$1 \times m$ and $1 \times k$ row vectors, respectively. It
follows that the post-multiplying matrix carries out the following
action:\\
(i) the columns involving $A_k$ remain unchanged; \\
(ii) column $1$ of $\left[ \begin{array}{ccc} \gamma & \beta &
\alpha \\ 0 & B_k & A_k \end{array} \right]$ having a pivot
in its first position, is moved into column $\mu(k+1)$;\\
(iii) column $\mu(k+1)+1$ of $\left[ \begin{array}{ccc} \gamma &
\beta & \alpha \\ 0 & B_k & A_k \end{array} \right]$ is
moved into column $m+1$;\\
(iv) columns $2,\ldots,\mu(k+1)$ and $\mu(k+1)+2,\ldots,m+1$ of
$\left[ \begin{array}{ccc} \gamma & \beta & \alpha \\ 0 & B_k &
A_k \end{array} \right]$ are moved one position to the left, into
columns $1,\ldots,\mu(k+1)-1$ and $\mu(k+1)+1,\ldots,m$,
respectively.
\\
This shows that the pivot structure of $B_{k+1}$ is obtained from
the pivot structure of $B_k$ in the following way: all the nonzero
entries of the structure vector $\eta_k$ are increased by $1$
except for the entry with index $\mu(k+1)$, which is reset to $1$.
This means that the pivot structure of $B_{k+1}$ is indeed given
by the vector $\eta_{k+1}$.

It remains to show that the matrix $A_{k+1}$ again has a staircase
form. Now, $A_k$ has a staircase form according to the induction
hypothesis and $A_{k+1}$ is recognized to be of the form
\[
  A_{k+1} = \left[ \begin{array}{cc}
  \delta & \alpha \\
  \epsilon & A_k
  \end{array} \right]
\]
where $\delta$ is the $\mu(k+1)$-st entry of the row vector
$\beta$ and $\epsilon$ is the $\mu(k+1)$-st column of $B_k$. This
means that the pivot in the first column of $A_{k+1}$ shows up in
position $(\eta_k)_{\mu(k+1)}+1$, which is equal to $\xi_k+1$.
However, it has already been established that the sequence
$\xi_{n-1},\xi_{n-2},\ldots,\xi_{p_B}$ is increasing. Therefore,
$A_{k+1}$ also has a staircase form.

By induction this shows for all $k=1,2,\ldots,n$, that $[B_k,A_k]$
has an admissible pivot structure for which the vector $\eta_k$
specifies the pivot structure of the matrix $B_k$. In particular,
for $k=n$ the claim of the theorem follows.
$\mbox{}$ \hfill $\Box$%

\section{An atlas for input-normal pairs $(A,B)$ under orthogonal
state-space equivalence, with $m=3$ and $n=4$}
\label{App:atlas}%

To illustrate the results and constructions of this paper, we here
present an atlas for the manifold of (controllable) input-normal
pairs $(A,B)$ under orthogonal state-space equivalence, for the
non-trivial case $m=3$ and $n=4$. Each of the charts in this atlas
gives rise to a particular full pivot structure in the
controllability matrix $K$ and an admissible pivot structure for
the row-orthonormal matrix $[B,A]$.

For given $m$ and $n$, the number of different \emph{admissible
numbered Young diagrams} (see the end of Section
\ref{triangularstructures}) is specified by
$\sum_{\ell=1}^{\min\{m,n\}} \ell ! \left( \begin{array}{c} m \\
\ell \end{array} \right) \left( \begin{array}{c} n-1 \\ \ell-1
\end{array} \right)$. For the case $m=3$ and $n=4$ this amounts to
$39$. To obtain a \emph{minimal} sub-atlas, precisely one chart
should be included for each nice selection, i.e. for each vector
of dynamical indices $d$ in ${\cal D}(m,n)$. The cardinality of
${\cal D}(m,n)$ is easily computed as $\left(
\begin{array}{c} m+n-1 \\ m-1 \end{array} \right)$. For the case
$m=3$ and $n=4$ this implies that a minimal sub-atlas consists of
$15$ charts. In Tables \ref{tab:part1}-\ref{tab:part3} the 15
different vectors of dynamical indices for this example are
displayed, along with the corresponding 39 admissible numbered
Young diagrams and their associated pivot structures in $K$ and in
$[B,A]$.

To arrive at an explicit \emph{parameterization} of a chart in
these tables, one may proceed in the \emph{discrete-time} case by
exploiting Eqn.\ (\ref{matrixproductformula}) for the construction
of orthogonal realization matrices, corresponding to balanced
realizations of discrete-time lossless systems. Here the sequence
of direction vectors $\{u_1,u_2,u_3,u_4\}$ is chosen to consist of
particular standard basis vectors, as indicated for each chart in
these tables too. The parameters are then provided by the sequence
of Schur vectors $\{v_1,v_2,v_3,v_4\}$ which are all required to
be of length $<1$. The $3 \times 3$ orthogonal matrix block $D_0$
can be set to any fixed value; the choice $D_0=I_3$ is a
convenient one. The latter is a consequence of the general fact
that if an orthogonal realization matrix $R = \left[
\begin{array}{cc} D & C \\ B & A \end{array} \right]$ is generated
by Eqn.\ (\ref{matrixproductformula}) for some $D_0$,
$\{u_1,\ldots,u_n\}$ and $\{v_1,\ldots,v_n\}$, then the
alternative choices $I_m$, $\{u_1,\ldots,u_n\}$ and $\{D_0^T
v_1,\ldots,D_0^T v_n\}$ yield the realization matrix $\left[
\begin{array}{cc} D_0 D & D_0 C \\ B & A \end{array} \right]$,
which exhibits exactly the \emph{same} input pair $(A,B)$.

From such a (minimal or non-minimal) atlas for input-normal pairs
under orthogonal state-space equivalence, a corresponding atlas
for all input-normal pairs of the given dimensions $m$ and $n$ is
directly obtained by regarding the associated manifold as a
Cartesian product of the previous manifold with the orthogonal
group ${\cal O}(n)$, related to the choice of state-space
transformation.

To arrive at a corresponding atlas for $m \times m$ lossless
systems of order $n$ one may instead regard this space as a
Cartesian product of the previous manifold with the orthogonal
group ${\cal O}(m)$, now related to the choice of $D_0$.

To obtain an atlas for asymptotically stable discrete-time systems
of order $n$ with $m$ inputs and $p$ outputs, one may proceed by
taking all the entries of $C$ and $D$ (of sizes $p \times n$ and
$p \times m$, respectively) to be free parameters, only subject to
the constraint that observability needs to hold for the pair
$(C,A)$ (a property which is then generically satisfied in each
chart, i.e. it only excludes a thin subset of parameter vectors).
Such an approach is useful in system identification, for instance
in conjunction with the method of separable least-squares (see
\cite{BCHV1999}). Then we may have to consider output-normal forms
instead, but this can be achieved easily using input-output
duality.

Finally, to deal with the \emph{continuous-time} case, the
well-known bilinear transform can of course be applied. However,
this will in general destroy the pivot structure in $K$ and in
$[B,A]$. To employ the results \emph{directly} in the
continuous-time case too, note that the pivot structures for
(controllable) input-normal pairs $[B,A]$ as given in the Tables
\ref{tab:part1}-\ref{tab:part3} do in fact apply to the
continuous-time case already, giving rise to local canonical forms
that can be computed numerically for a given state-space
realization in a straightforward way. What at present seems to be
lacking in the continuous-time case is an explicit
parameterization of these local canonical forms (such as may be
required in system identification). This is currently the topic of
ongoing research.

\begin{table}[htbp]
\tiny
    \centering
        \begin{tabular}{ccccccc}
      {C}hart & Young diagram $Z$ and & Admissible numbered & Full pivot                        &
      Sequence of direction         & Admissible pivot               & Structure of \\
              & dynamical index $d$   & Young diagrams $Y$  & structure $\widetilde{J}$ for $K$ &
      vectors $\{u_1,u_2,u_3,u_4\}$ & structure $J$ for $[B,A]$      & $[B,A]$      \\ \hline \hline
& & & & & & \\
   1 & $\begin{array}{|c|c|c|c|} \hline
         1 & 1 & 1 & 1 \\ \hline
           &   &   &   \\ \hline
           &   &   &   \\ \hline
       \end{array}$ &
       $\begin{array}{|c|c|c|c|} \hline
         1 & 2 & 3 & 4 \\ \hline
           &   &   &   \\ \hline
           &   &   &   \\ \hline
       \end{array}$ &
       $\{1,4,7,10\}$ &
       $\{e_1,e_1,e_1,e_1\}$ &
       $\{1,4,5,6\}$ &
       $\left[ \begin{array}{ccc|cccc}
         {+} & \ast & \ast & \ast & \ast & \ast & \ast \\
          0  & \ast & \ast &  +   & \ast & \ast & \ast \\
          0  & \ast & \ast &  0   &  +   & \ast & \ast \\
          0  & \ast & \ast &  0   &  0   &  +   & \ast \\
       \end{array} \right]$ \\
& & & & & & \\
   & $(4,0,0)$ & & & & & \\
& & & & & & \\ \hline
& & & & & & \\
 2 & $\begin{array}{|c|c|c|c|} \hline
         1 & 1 & 1 & \phantom{1} \\ \hline
         1 &   &   &             \\ \hline
           &   &   &             \\ \hline
       \end{array}$ &
       $\begin{array}{|c|c|c|c|} \hline
         1 & 2 & 3 & \phantom{4} \\ \hline
         4 &   &   &   \\ \hline
           &   &   &   \\ \hline
       \end{array}$ &
       $\{1,4,7,2\}$ &
       $\{e_2,e_1,e_1,e_1\}$ &
       $\{1,4,5,2\}$ &
       $\left[ \begin{array}{ccc|cccc}
         {+} & \ast & \ast & \ast & \ast & \ast & \ast \\
          0  & \ast & \ast &  +   & \ast & \ast & \ast \\
          0  & \ast & \ast &  0   &  +   & \ast & \ast \\
          0  &  +   & \ast &  0   &  0   & \ast & \ast \\
       \end{array} \right]$ \\
& & & & & & \\
   & $(3,1,0)$ & $\begin{array}{|c|c|c|c|} \hline
         1 & 2 & 4 & \phantom{4} \\ \hline
         3 &   &   &   \\ \hline
           &   &   &   \\ \hline
       \end{array}$ &
       $\{1,4,2,7\}$ &
       $\{e_1,e_2,e_1,e_1\}$ &
       $\{1,4,2,5\}$ &
       $\left[ \begin{array}{ccc|cccc}
         {+} & \ast & \ast & \ast & \ast & \ast & \ast \\
          0  & \ast & \ast &  +   & \ast & \ast & \ast \\
          0  &  +   & \ast &  0   & \ast & \ast & \ast \\
          0  &  0   & \ast &  0   &  +   & \ast & \ast \\
       \end{array} \right]$ \\
& & & & & & \\ \hline
& & & & & & \\
 3 & $\begin{array}{|c|c|c|c|} \hline
         1 & 1 & 1 & \phantom{1} \\ \hline
           &   &   &             \\ \hline
         1 &   &   &             \\ \hline
       \end{array}$ &
       $\begin{array}{|c|c|c|c|} \hline
         1 & 2 & 3 & \phantom{4} \\ \hline
           &   &   &   \\ \hline
         4 &   &   &   \\ \hline
       \end{array}$ &
       $\{1,4,7,3\}$ &
       $\{e_3,e_1,e_1,e_1\}$ &
       $\{1,4,5,3\}$ &
       $\left[ \begin{array}{ccc|cccc}
         {+} & \ast & \ast & \ast & \ast & \ast & \ast \\
          0  & \ast & \ast &  +   & \ast & \ast & \ast \\
          0  & \ast & \ast &  0   &  +   & \ast & \ast \\
          0  & \ast &  +   &  0   &  0   & \ast & \ast \\
       \end{array} \right]$ \\
& & & & & & \\
   & $(3,0,1)$ & $\begin{array}{|c|c|c|c|} \hline
         1 & 2 & 4 & \phantom{4} \\ \hline
           &   &   &   \\ \hline
         3 &   &   &   \\ \hline
       \end{array}$ &
       $\{1,4,3,7\}$ &
       $\{e_1,e_3,e_1,e_1\}$ &
       $\{1,4,3,5\}$ &
       $\left[ \begin{array}{ccc|cccc}
         {+} & \ast & \ast & \ast & \ast & \ast & \ast \\
          0  & \ast & \ast &  +   & \ast & \ast & \ast \\
          0  & \ast &  +   &  0   & \ast & \ast & \ast \\
          0  & \ast &  0   &  0   &  +   & \ast & \ast \\
       \end{array} \right]$ \\
& & & & & & \\ \hline
& & & & & & \\
 4 & $\begin{array}{|c|c|c|c|} \hline
         1 & 1 & \phantom{1} & \phantom{1} \\ \hline
         1 & 1 &             &             \\ \hline
           &   &             &             \\ \hline
       \end{array}$ &
       $\begin{array}{|c|c|c|c|} \hline
         1 & 3 & \phantom{3} & \phantom{4} \\ \hline
         2 & 4 &   &   \\ \hline
           &   &   &   \\ \hline
       \end{array}$ &
       $\{1,2,4,5\}$ &
       $\{e_2,e_1,e_2,e_1\}$ &
       $\{1,2,4,5\}$ &
       $\left[ \begin{array}{ccc|cccc}
         {+} & \ast & \ast & \ast & \ast & \ast & \ast \\
          0  &  +   & \ast & \ast & \ast & \ast & \ast \\
          0  &  0   & \ast &  +   & \ast & \ast & \ast \\
          0  &  0   & \ast &  0   &  +   & \ast & \ast \\
       \end{array} \right]$ \\
& & & & & & \\
   & $(2,2,0)$ & $\begin{array}{|c|c|c|c|} \hline
         2 & 4 & \phantom{3} & \phantom{4} \\ \hline
         1 & 3 &   &   \\ \hline
           &   &   &   \\ \hline
       \end{array}$ &
       $\{2,1,5,4\}$ &
       $\{e_1,e_2,e_1,e_2\}$ &
       $\{2,1,4,5\}$ &
       $\left[ \begin{array}{ccc|cccc}
        \ast &  +   & \ast & \ast & \ast & \ast & \ast \\
         {+} &  0   & \ast & \ast & \ast & \ast & \ast \\
          0  &  0   & \ast &  +   & \ast & \ast & \ast \\
          0  &  0   & \ast &  0   &  +   & \ast & \ast \\
       \end{array} \right]$ \\
& & & & & & \\ \hline
& & & & & & \\
 5 & $\begin{array}{|c|c|c|c|} \hline
         1 & 1 & \phantom{1} & \phantom{1} \\ \hline
         1 &   &             &             \\ \hline
         1 &   &             &             \\ \hline
       \end{array}$ &
       $\begin{array}{|c|c|c|c|} \hline
         1 & 2 & \phantom{3} & \phantom{4} \\ \hline
         3 &   &   &   \\ \hline
         4 &   &   &   \\ \hline
       \end{array}$ &
       $\{1,4,2,3\}$ &
       $\{e_3,e_2,e_1,e_1\}$ &
       $\{1,4,2,3\}$ &
       $\left[ \begin{array}{ccc|cccc}
         {+} & \ast & \ast & \ast & \ast & \ast & \ast \\
          0  & \ast & \ast &  +   & \ast & \ast & \ast \\
          0  &  +   & \ast &  0   & \ast & \ast & \ast \\
          0  &  0   &  +   &  0   & \ast & \ast & \ast \\
       \end{array} \right]$ \\
& & & & & & \\
   & $(2,1,1)$ & $\begin{array}{|c|c|c|c|} \hline
         1 & 2 & \phantom{3} & \phantom{4} \\ \hline
         4 &   &   &   \\ \hline
         3 &   &   &   \\ \hline
       \end{array}$ &
       $\{1,4,3,2\}$ &
       $\{e_2,e_3,e_1,e_1\}$ &
       $\{1,4,3,2\}$ &
       $\left[ \begin{array}{ccc|cccc}
         {+} & \ast & \ast & \ast & \ast & \ast & \ast \\
          0  & \ast & \ast &  +   & \ast & \ast & \ast \\
          0  & \ast &  +   &  0   & \ast & \ast & \ast \\
          0  &  +   &  0   &  0   & \ast & \ast & \ast \\
       \end{array} \right]$ \\
& & & & & & \\
   & & $\begin{array}{|c|c|c|c|} \hline
         1 & 3 & \phantom{3} & \phantom{4} \\ \hline
         2 &   &   &   \\ \hline
         4 &   &   &   \\ \hline
       \end{array}$ &
       $\{1,2,4,3\}$ &
       $\{e_3,e_1,e_2,e_1\}$ &
       $\{1,2,4,3\}$ &
       $\left[ \begin{array}{ccc|cccc}
         {+} & \ast & \ast & \ast & \ast & \ast & \ast \\
          0  &  +   & \ast & \ast & \ast & \ast & \ast \\
          0  &  0   & \ast &  +   & \ast & \ast & \ast \\
          0  &  0   &  +   &  0   & \ast & \ast & \ast \\
       \end{array} \right]$ \\
& & & & & & \\
   & & $\begin{array}{|c|c|c|c|} \hline
         1 & 3 & \phantom{3} & \phantom{4} \\ \hline
         4 &   &   &   \\ \hline
         2 &   &   &   \\ \hline
       \end{array}$ &
       $\{1,3,4,2\}$ &
       $\{e_2,e_1,e_3,e_1\}$ &
       $\{1,3,4,2\}$ &
       $\left[ \begin{array}{ccc|cccc}
         {+} & \ast & \ast & \ast & \ast & \ast & \ast \\
          0  & \ast &  +   & \ast & \ast & \ast & \ast \\
          0  & \ast &  0   &  +   & \ast & \ast & \ast \\
          0  &  +   &  0   &  0   & \ast & \ast & \ast \\
       \end{array} \right]$ \\
& & & & & & \\
   & & $\begin{array}{|c|c|c|c|} \hline
         1 & 4 & \phantom{3} & \phantom{4} \\ \hline
         2 &   &   &   \\ \hline
         3 &   &   &   \\ \hline
       \end{array}$ &
       $\{1,2,3,4\}$ &
       $\{e_1,e_3,e_2,e_1\}$ &
       $\{1,2,3,4\}$ &
       $\left[ \begin{array}{ccc|cccc}
         {+} & \ast & \ast & \ast & \ast & \ast & \ast \\
          0  &  +   & \ast & \ast & \ast & \ast & \ast \\
          0  &  0   &  +   & \ast & \ast & \ast & \ast \\
          0  &  0   &  0   &  +   & \ast & \ast & \ast \\
       \end{array} \right]$ \\
& & & & & & \\
   & & $\begin{array}{|c|c|c|c|} \hline
         1 & 4 & \phantom{3} & \phantom{4} \\ \hline
         3 &   &   &   \\ \hline
         2 &   &   &   \\ \hline
       \end{array}$ &
       $\{1,3,2,4\}$ &
       $\{e_1,e_2,e_3,e_1\}$ &
       $\{1,3,2,4\}$ &
       $\left[ \begin{array}{ccc|cccc}
         {+} & \ast & \ast & \ast & \ast & \ast & \ast \\
          0  & \ast &  +   & \ast & \ast & \ast & \ast \\
          0  &  +   &  0   & \ast & \ast & \ast & \ast \\
          0  &  0   &  0   &  +   & \ast & \ast & \ast \\
       \end{array} \right]$ \\
       & & & & & & \\      \hline
        \end{tabular}
    \caption{Charts 1--5 for $[B,A]$, for the case $m=3$ and $n=4$.}
    \label{tab:part1}
\end{table}

\begin{table}[htbp]
\tiny
    \centering
        \begin{tabular}{ccccccc}
      {C}hart & Young diagram $Z$ and & Admissible numbered & Full pivot                        &
      Sequence of direction         & Admissible pivot               & Structure of \\
              & dynamical index $d$   & Young diagrams $Y$  & structure $\widetilde{J}$ for $K$ &
      vectors $\{u_1,u_2,u_3,u_4\}$ & structure $J$ for $[B,A]$      & $[B,A]$      \\ \hline \hline
& & & & & & \\
   6 & $\begin{array}{|c|c|c|c|} \hline
           &   &   &   \\ \hline
         1 & 1 & 1 & 1 \\ \hline
           &   &   &   \\ \hline
       \end{array}$ &
       $\begin{array}{|c|c|c|c|} \hline
           &   &   &   \\ \hline
         1 & 2 & 3 & 4 \\ \hline
           &   &   &   \\ \hline
       \end{array}$ &
       $\{2,5,8,11\}$ &
       $\{e_2,e_2,e_2,e_2\}$ &
       $\{2,4,5,6\}$ &
       $\left[ \begin{array}{ccc|cccc}
        \ast &  +   & \ast & \ast & \ast & \ast & \ast \\
        \ast &  0   & \ast &  +   & \ast & \ast & \ast \\
        \ast &  0   & \ast &  0   &  +   & \ast & \ast \\
        \ast &  0   & \ast &  0   &  0   &  +   & \ast \\
       \end{array} \right]$ \\
& & & & & & \\
   & $(0,4,0)$ & & & & & \\
& & & & & & \\ \hline
& & & & & & \\
 7 & $\begin{array}{|c|c|c|c|} \hline
           &   &   &             \\ \hline
         1 & 1 & 1 & \phantom{1} \\ \hline
         1 &   &   &             \\ \hline
       \end{array}$ &
       $\begin{array}{|c|c|c|c|} \hline
           &   &   &   \\ \hline
         1 & 2 & 3 & \phantom{4} \\ \hline
         4 &   &   &   \\ \hline
       \end{array}$ &
       $\{2,5,8,3\}$ &
       $\{e_3,e_2,e_2,e_2\}$ &
       $\{2,4,5,3\}$ &
       $\left[ \begin{array}{ccc|cccc}
         \ast &  +  & \ast & \ast & \ast & \ast & \ast \\
         \ast &  0  & \ast &  +   & \ast & \ast & \ast \\
         \ast &  0  & \ast &  0   &  +   & \ast & \ast \\
         \ast &  0  &  +   &  0   &  0   & \ast & \ast \\
       \end{array} \right]$ \\
& & & & & & \\
   & $(0,3,1)$ & $\begin{array}{|c|c|c|c|} \hline
           &   &   &   \\ \hline
         1 & 2 & 4 & \phantom{4} \\ \hline
         3 &   &   &   \\ \hline
       \end{array}$ &
       $\{2,5,3,8\}$ &
       $\{e_2,e_3,e_2,e_2\}$ &
       $\{2,4,3,5\}$ &
       $\left[ \begin{array}{ccc|cccc}
         \ast &  +  & \ast & \ast & \ast & \ast & \ast \\
         \ast &  0  & \ast &  +   & \ast & \ast & \ast \\
         \ast &  0  &  +   &  0   & \ast & \ast & \ast \\
         \ast &  0  &  0   &  0   &  +   & \ast & \ast \\
       \end{array} \right]$ \\
& & & & & & \\ \hline
& & & & & & \\
 8 & $\begin{array}{|c|c|c|c|} \hline
         1 &   &   &             \\ \hline
         1 & 1 & 1 & \phantom{1} \\ \hline
           &   &   &             \\ \hline
       \end{array}$ &
       $\begin{array}{|c|c|c|c|} \hline
         4 &   &   &   \\ \hline
         1 & 2 & 3 & \phantom{4} \\ \hline
           &   &   &   \\ \hline
       \end{array}$ &
       $\{2,5,8,1\}$ &
       $\{e_1,e_2,e_2,e_2\}$ &
       $\{2,4,5,1\}$ &
       $\left[ \begin{array}{ccc|cccc}
         \ast &  +  & \ast & \ast & \ast & \ast & \ast \\
         \ast &  0  & \ast &  +   & \ast & \ast & \ast \\
         \ast &  0  & \ast &  0   &  +   & \ast & \ast \\
         {+}  &  0  & \ast &  0   &  0   & \ast & \ast \\
       \end{array} \right]$ \\
& & & & & & \\
   & $(1,3,0)$ & $\begin{array}{|c|c|c|c|} \hline
         3 &   &   &   \\ \hline
         1 & 2 & 4 & \phantom{4} \\ \hline
           &   &   &   \\ \hline
       \end{array}$ &
       $\{2,5,1,8\}$ &
       $\{e_2,e_1,e_2,e_2\}$ &
       $\{2,4,1,5\}$ &
       $\left[ \begin{array}{ccc|cccc}
         \ast &  +  & \ast & \ast & \ast & \ast & \ast \\
         \ast &  0  & \ast &  +   & \ast & \ast & \ast \\
         {+}  &  0  & \ast &  0   & \ast & \ast & \ast \\
          0   &  0  & \ast &  0   &  +   & \ast & \ast \\
       \end{array} \right]$ \\
& & & & & & \\ \hline
& & & & & & \\
 9 & $\begin{array}{|c|c|c|c|} \hline
           &   &             &             \\ \hline
         1 & 1 & \phantom{1} & \phantom{1} \\ \hline
         1 & 1 &             &             \\ \hline
       \end{array}$ &
       $\begin{array}{|c|c|c|c|} \hline
           &   &   &   \\ \hline
         1 & 3 & \phantom{3} & \phantom{4} \\ \hline
         2 & 4 &   &   \\ \hline
       \end{array}$ &
       $\{2,3,5,6\}$ &
       $\{e_3,e_2,e_3,e_2\}$ &
       $\{2,3,4,5\}$ &
       $\left[ \begin{array}{ccc|cccc}
         \ast &  +  & \ast & \ast & \ast & \ast & \ast \\
         \ast &  0  &  +   & \ast & \ast & \ast & \ast \\
         \ast &  0  &  0   &  +   & \ast & \ast & \ast \\
         \ast &  0  &  0   &  0   &  +   & \ast & \ast \\
       \end{array} \right]$ \\
& & & & & & \\
   & $(0,2,2)$ & $\begin{array}{|c|c|c|c|} \hline
           &   &   &   \\ \hline
         2 & 4 & \phantom{3} & \phantom{4} \\ \hline
         1 & 3 &   &   \\ \hline
       \end{array}$ &
       $\{3,2,6,5\}$ &
       $\{e_2,e_3,e_2,e_3\}$ &
       $\{3,2,4,5\}$ &
       $\left[ \begin{array}{ccc|cccc}
        \ast & \ast &  +   & \ast & \ast & \ast & \ast \\
        \ast &   +  &  0   & \ast & \ast & \ast & \ast \\
        \ast &   0  &  0   &  +   & \ast & \ast & \ast \\
        \ast &   0  &  0   &  0   &  +   & \ast & \ast \\
       \end{array} \right]$ \\
& & & & & & \\ \hline
& & & & & & \\
 10 & $\begin{array}{|c|c|c|c|} \hline
         1 &   &             &             \\ \hline
         1 & 1 & \phantom{1} & \phantom{1} \\ \hline
         1 &   &             &             \\ \hline
       \end{array}$ &
       $\begin{array}{|c|c|c|c|} \hline
         4 &   &   &   \\ \hline
         1 & 2 & \phantom{3} & \phantom{4} \\ \hline
         3 &   &   &   \\ \hline
       \end{array}$ &
       $\{2,5,3,1\}$ &
       $\{e_1,e_3,e_2,e_2\}$ &
       $\{2,4,3,1\}$ &
       $\left[ \begin{array}{ccc|cccc}
         \ast &  +  & \ast & \ast & \ast & \ast & \ast \\
         \ast &  0  & \ast &  +   & \ast & \ast & \ast \\
         \ast &  0  &  +   &  0   & \ast & \ast & \ast \\
         {+}  &  0  &  0   &  0   & \ast & \ast & \ast \\
       \end{array} \right]$ \\
& & & & & & \\
   & $(1,2,1)$ & $\begin{array}{|c|c|c|c|} \hline
         3 &   &   &   \\ \hline
         1 & 2 & \phantom{3} & \phantom{4} \\ \hline
         4 &   &   &   \\ \hline
       \end{array}$ &
       $\{2,5,1,3\}$ &
       $\{e_3,e_1,e_2,e_2\}$ &
       $\{2,4,1,3\}$ &
       $\left[ \begin{array}{ccc|cccc}
         \ast &  +  & \ast & \ast & \ast & \ast & \ast \\
         \ast &  0  & \ast &  +   & \ast & \ast & \ast \\
         {+}  &  0  & \ast &  0   & \ast & \ast & \ast \\
          0   &  0  &  +   &  0   & \ast & \ast & \ast \\
       \end{array} \right]$ \\
& & & & & & \\
   & & $\begin{array}{|c|c|c|c|} \hline
         4 &   &   &   \\ \hline
         1 & 3 & \phantom{3} & \phantom{4} \\ \hline
         2 &   &   &   \\ \hline
       \end{array}$ &
       $\{2,3,5,1\}$ &
       $\{e_1,e_2,e_3,e_2\}$ &
       $\{2,3,4,1\}$ &
       $\left[ \begin{array}{ccc|cccc}
         \ast &  +  & \ast & \ast & \ast & \ast & \ast \\
         \ast &  0  &  +   & \ast & \ast & \ast & \ast \\
         \ast &  0  &  0   &  +   & \ast & \ast & \ast \\
         {+}  &  0  &  0   &  0   & \ast & \ast & \ast \\
       \end{array} \right]$ \\
& & & & & & \\
   & & $\begin{array}{|c|c|c|c|} \hline
         2 &   &   &   \\ \hline
         1 & 3 & \phantom{3} & \phantom{4} \\ \hline
         4 &   &   &   \\ \hline
       \end{array}$ &
       $\{2,1,5,3\}$ &
       $\{e_3,e_2,e_1,e_2\}$ &
       $\{2,1,4,3\}$ &
       $\left[ \begin{array}{ccc|cccc}
         \ast &  +  & \ast & \ast & \ast & \ast & \ast \\
         {+}  &  0  & \ast & \ast & \ast & \ast & \ast \\
          0   &  0  & \ast &  +   & \ast & \ast & \ast \\
          0   &  0  &  +   &  0   & \ast & \ast & \ast \\
       \end{array} \right]$ \\
& & & & & & \\
   & & $\begin{array}{|c|c|c|c|} \hline
         3 &   &   &   \\ \hline
         1 & 4 & \phantom{3} & \phantom{4} \\ \hline
         2 &   &   &   \\ \hline
       \end{array}$ &
       $\{2,3,1,5\}$ &
       $\{e_2,e_1,e_3,e_2\}$ &
       $\{2,3,1,4\}$ &
       $\left[ \begin{array}{ccc|cccc}
         \ast &  +  & \ast & \ast & \ast & \ast & \ast \\
         \ast &  0  &  +   & \ast & \ast & \ast & \ast \\
         {+}  &  0  &  0   & \ast & \ast & \ast & \ast \\
          0   &  0  &  0   &  +   & \ast & \ast & \ast \\
       \end{array} \right]$ \\
& & & & & & \\
   & & $\begin{array}{|c|c|c|c|} \hline
         2 &   &   &   \\ \hline
         1 & 4 & \phantom{3} & \phantom{4} \\ \hline
         3 &   &   &   \\ \hline
       \end{array}$ &
       $\{2,1,3,5\}$ &
       $\{e_2,e_3,e_1,e_2\}$ &
       $\{2,1,3,4\}$ &
       $\left[ \begin{array}{ccc|cccc}
         \ast &  +  & \ast & \ast & \ast & \ast & \ast \\
         {+}  &  0  & \ast & \ast & \ast & \ast & \ast \\
          0   &  0  &  +   & \ast & \ast & \ast & \ast \\
          0   &  0  &  0   &  +   & \ast & \ast & \ast \\
       \end{array} \right]$ \\
& & & & & & \\ \hline
        \end{tabular}
    \caption{Charts 6--10 for $[B,A]$, for the case $m=3$ and $n=4$.}
    \label{tab:part2}
\end{table}

\begin{table}[htbp]
\tiny
    \centering
        \begin{tabular}{ccccccc}
      {C}hart & Young diagram $Z$ and & Admissible numbered & Full pivot                        &
      Sequence of direction         & Admissible pivot               & Structure of \\
              & dynamical index $d$   & Young diagrams $Y$  & structure $\widetilde{J}$ for $K$ &
      vectors $\{u_1,u_2,u_3,u_4\}$ & structure $J$ for $[B,A]$      & $[B,A]$      \\ \hline \hline
& & & & & & \\
   11 & $\begin{array}{|c|c|c|c|} \hline
           &   &   &   \\ \hline
           &   &   &   \\ \hline
         1 & 1 & 1 & 1 \\ \hline
       \end{array}$ &
       $\begin{array}{|c|c|c|c|} \hline
           &   &   &   \\ \hline
           &   &   &   \\ \hline
         1 & 2 & 3 & 4 \\ \hline
       \end{array}$ &
       $\{3,6,9,12\}$ &
       $\{e_3,e_3,e_3,e_3\}$ &
       $\{3,4,5,6\}$ &
       $\left[ \begin{array}{ccc|cccc}
         \ast & \ast &  +  & \ast & \ast & \ast & \ast \\
         \ast & \ast &  0  &  +   & \ast & \ast & \ast \\
         \ast & \ast &  0  &  0   &  +   & \ast & \ast \\
         \ast & \ast &  0  &  0   &  0   &  +   & \ast \\
       \end{array} \right]$ \\
& & & & & & \\
   & $(0,0,4)$ & & & & & \\
& & & & & & \\ \hline
& & & & & & \\
 12 & $\begin{array}{|c|c|c|c|} \hline
         1 &   &   &             \\ \hline
           &   &   &             \\ \hline
         1 & 1 & 1 & \phantom{1} \\ \hline
       \end{array}$ &
       $\begin{array}{|c|c|c|c|} \hline
         4 &   &   &   \\ \hline
           &   &   &   \\ \hline
         1 & 2 & 3 & \phantom{4} \\ \hline
       \end{array}$ &
       $\{3,6,9,1\}$ &
       $\{e_1,e_3,e_3,e_3\}$ &
       $\{3,4,5,1\}$ &
       $\left[ \begin{array}{ccc|cccc}
         \ast & \ast &  +  & \ast & \ast & \ast & \ast \\
         \ast & \ast &  0  &  +   & \ast & \ast & \ast \\
         \ast & \ast &  0  &  0   &  +   & \ast & \ast \\
         {+}  & \ast &  0  &  0   &  0   & \ast & \ast \\
       \end{array} \right]$ \\
& & & & & & \\
   & $(1,0,3)$ & $\begin{array}{|c|c|c|c|} \hline
         3 &   &   &   \\ \hline
           &   &   &   \\ \hline
         1 & 2 & 4 & \phantom{4} \\ \hline
       \end{array}$ &
       $\{3,6,1,9\}$ &
       $\{e_3,e_1,e_3,e_3\}$ &
       $\{3,4,1,5\}$ &
       $\left[ \begin{array}{ccc|cccc}
         \ast & \ast &  +  & \ast & \ast & \ast & \ast \\
         \ast & \ast &  0  &  +   & \ast & \ast & \ast \\
         {+}  & \ast &  0  &  0   & \ast & \ast & \ast \\
          0   & \ast &  0  &  0   &  +   & \ast & \ast \\
       \end{array} \right]$ \\
& & & & & & \\ \hline
& & & & & & \\
 13 & $\begin{array}{|c|c|c|c|} \hline
           &   &   &             \\ \hline
         1 &   &   &             \\ \hline
         1 & 1 & 1 & \phantom{1} \\ \hline
       \end{array}$ &
       $\begin{array}{|c|c|c|c|} \hline
           &   &   &   \\ \hline
         4 &   &   &   \\ \hline
         1 & 2 & 3 & \phantom{4} \\ \hline
       \end{array}$ &
       $\{3,6,9,2\}$ &
       $\{e_2,e_3,e_3,e_3\}$ &
       $\{3,4,5,2\}$ &
       $\left[ \begin{array}{ccc|cccc}
         \ast & \ast &  +  & \ast & \ast & \ast & \ast \\
         \ast & \ast &  0  &  +   & \ast & \ast & \ast \\
         \ast & \ast &  0  &  0   &  +   & \ast & \ast \\
         \ast &  +   &  0  &  0   &  0   & \ast & \ast \\
       \end{array} \right]$ \\
& & & & & & \\
   & $(0,1,3)$ & $\begin{array}{|c|c|c|c|} \hline
           &   &   &   \\ \hline
         3 &   &   &   \\ \hline
         1 & 2 & 4 & \phantom{4} \\ \hline
       \end{array}$ &
       $\{3,6,2,9\}$ &
       $\{e_3,e_2,e_3,e_3\}$ &
       $\{3,4,2,5\}$ &
       $\left[ \begin{array}{ccc|cccc}
         \ast & \ast &  +  & \ast & \ast & \ast & \ast \\
         \ast & \ast &  0  &  +   & \ast & \ast & \ast \\
         \ast &  +   &  0  &  0   & \ast & \ast & \ast \\
         \ast &  0   &  0  &  0   &  +   & \ast & \ast \\
       \end{array} \right]$ \\
& & & & & & \\ \hline
& & & & & & \\
 14 & $\begin{array}{|c|c|c|c|} \hline
         1 & 1 &             &             \\ \hline
           &   &             &             \\ \hline
         1 & 1 & \phantom{1} & \phantom{1} \\ \hline
       \end{array}$ &
       $\begin{array}{|c|c|c|c|} \hline
         2 & 4 &   &   \\ \hline
           &   &   &   \\ \hline
         1 & 3 & \phantom{3} & \phantom{4} \\ \hline
       \end{array}$ &
       $\{3,1,6,4\}$ &
       $\{e_1,e_3,e_1,e_3\}$ &
       $\{3,1,4,5\}$ &
       $\left[ \begin{array}{ccc|cccc}
         \ast & \ast &  +  & \ast & \ast & \ast & \ast \\
         {+}  & \ast &  0  & \ast & \ast & \ast & \ast \\
          0   & \ast &  0  &  +   & \ast & \ast & \ast \\
          0   & \ast &  0  &  0   &  +   & \ast & \ast \\
       \end{array} \right]$ \\
& & & & & & \\
   & $(2,0,2)$ & $\begin{array}{|c|c|c|c|} \hline
         1 & 3 &   &   \\ \hline
           &   &   &   \\ \hline
         2 & 4 & \phantom{3} & \phantom{4} \\ \hline
       \end{array}$ &
       $\{1,3,4,6\}$ &
       $\{e_3,e_1,e_3,e_1\}$ &
       $\{1,3,4,5\}$ &
       $\left[ \begin{array}{ccc|cccc}
         {+} & \ast & \ast & \ast & \ast & \ast & \ast \\
          0  & \ast &  +   & \ast & \ast & \ast & \ast \\
          0  & \ast &  0   &  +   & \ast & \ast & \ast \\
          0  & \ast &  0   &  0   &  +   & \ast & \ast \\
       \end{array} \right]$ \\
& & & & & & \\ \hline
& & & & & & \\
 15 & $\begin{array}{|c|c|c|c|} \hline
         1 &   &             &             \\ \hline
         1 &   &             &             \\ \hline
         1 & 1 & \phantom{1} & \phantom{1} \\ \hline
       \end{array}$ &
       $\begin{array}{|c|c|c|c|} \hline
         3 &   &   &   \\ \hline
         4 &   &   &   \\ \hline
         1 & 2 & \phantom{3} & \phantom{4} \\ \hline
       \end{array}$ &
       $\{3,6,1,2\}$ &
       $\{e_2,e_1,e_3,e_3\}$ &
       $\{3,4,1,2\}$ &
       $\left[ \begin{array}{ccc|cccc}
         \ast & \ast &  +  & \ast & \ast & \ast & \ast \\
         \ast & \ast &  0  &  +   & \ast & \ast & \ast \\
         {+}  & \ast &  0  &  0   & \ast & \ast & \ast \\
          0   &  +   &  0  &  0   & \ast & \ast & \ast \\
       \end{array} \right]$ \\
& & & & & & \\
   & $(1,1,2)$ & $\begin{array}{|c|c|c|c|} \hline
         4 &   &   &   \\ \hline
         3 &   &   &   \\ \hline
         1 & 2 & \phantom{3} & \phantom{4} \\ \hline
       \end{array}$ &
       $\{3,6,2,1\}$ &
       $\{e_1,e_2,e_3,e_3\}$ &
       $\{3,4,2,1\}$ &
       $\left[ \begin{array}{ccc|cccc}
         \ast & \ast &  +  & \ast & \ast & \ast & \ast \\
         \ast & \ast &  0  &  +   & \ast & \ast & \ast \\
         \ast &  +   &  0  &  0   & \ast & \ast & \ast \\
         {+}  &  0   &  0  &  0   & \ast & \ast & \ast \\
       \end{array} \right]$ \\
& & & & & & \\
   & & $\begin{array}{|c|c|c|c|} \hline
         2 &   &   &   \\ \hline
         4 &   &   &   \\ \hline
         1 & 3 & \phantom{3} & \phantom{4} \\ \hline
       \end{array}$ &
       $\{3,1,6,2\}$ &
       $\{e_2,e_3,e_1,e_3\}$ &
       $\{3,1,4,2\}$ &
       $\left[ \begin{array}{ccc|cccc}
         \ast & \ast &  +  & \ast & \ast & \ast & \ast \\
         {+}  & \ast &  0  & \ast & \ast & \ast & \ast \\
          0   & \ast &  0  &  +   & \ast & \ast & \ast \\
          0   &  +   &  0  &  0   & \ast & \ast & \ast \\
       \end{array} \right]$ \\
& & & & & & \\
   & & $\begin{array}{|c|c|c|c|} \hline
         4 &   &   &   \\ \hline
         2 &   &   &   \\ \hline
         1 & 3 & \phantom{3} & \phantom{4} \\ \hline
       \end{array}$ &
       $\{3,2,6,1\}$ &
       $\{e_1,e_3,e_2,e_3\}$ &
       $\{3,2,4,1\}$ &
       $\left[ \begin{array}{ccc|cccc}
         \ast & \ast &  +  & \ast & \ast & \ast & \ast \\
         \ast &  +   &  0  & \ast & \ast & \ast & \ast \\
         \ast &  0   &  0  &  +   & \ast & \ast & \ast \\
         {+}  &  0   &  0  &  0   & \ast & \ast & \ast \\
       \end{array} \right]$ \\
& & & & & & \\
   & & $\begin{array}{|c|c|c|c|} \hline
         2 &   &   &   \\ \hline
         3 &   &   &   \\ \hline
         1 & 4 & \phantom{3} & \phantom{4} \\ \hline
       \end{array}$ &
       $\{3,1,2,6\}$ &
       $\{e_3,e_2,e_1,e_3\}$ &
       $\{3,1,2,4\}$ &
       $\left[ \begin{array}{ccc|cccc}
         \ast & \ast &  +  & \ast & \ast & \ast & \ast \\
         {+}  & \ast &  0  & \ast & \ast & \ast & \ast \\
          0   &  +   &  0  & \ast & \ast & \ast & \ast \\
          0   &  0   &  0  &  +   & \ast & \ast & \ast \\
       \end{array} \right]$ \\
& & & & & & \\
   & & $\begin{array}{|c|c|c|c|} \hline
         3 &   &   &   \\ \hline
         2 &   &   &   \\ \hline
         1 & 4 & \phantom{3} & \phantom{4} \\ \hline
       \end{array}$ &
       $\{3,2,1,6\}$ &
       $\{e_3,e_1,e_2,e_3\}$ &
       $\{3,2,1,4\}$ &
       $\left[ \begin{array}{ccc|cccc}
         \ast & \ast &  +  & \ast & \ast & \ast & \ast \\
         \ast &  +   &  0  & \ast & \ast & \ast & \ast \\
         {+}  &  0   &  0  & \ast & \ast & \ast & \ast \\
          0   &  0   &  0  &  +   & \ast & \ast & \ast \\
       \end{array} \right]$ \\
& & & & & & \\ \hline
        \end{tabular}
    \caption{Charts 11--15 for $[B,A]$, for the case $m=3$ and $n=4$.}
    \label{tab:part3}
\end{table}

\end{document}